\documentclass[12pt,sort&compress]{elsarticle}

\graphicspath{{"Figures/"}}
\usepackage{amsmath,amssymb,enumerate,graphicx,subfig,float,color,mathtools,amsthm,booktabs,multirow,array,xcolor,colortbl,stmaryrd}
\usepackage{titlesec,bibentry}
\usepackage[section]{placeins} 
\usepackage{hyperref}
\newcommand{\revision}[1]{{#1}}
\newcommand{\change}[1]{{#1}}
\newcommand{\amend}[1]{{#1}}
\newcommand{\edit}[1]{{#1}}
\newcommand{\norm}[1]{\left\lVert#1\right\rVert}
\newcommand{\final}[1]{{#1}}
\newcommand{\finalrevision}[1]{{#1}}
\newcommand{\newedit}[1]{{#1}}
\newcommand{\neweredit}[1]{{#1}}
\newcommand{\finaledit}[1]{{#1}}
\newcommand{\finalfinaledit}[1]{{#1}}
\newcommand{\absolutefinaledit}[1]{{#1}}

\makeatletter
\def\ps@pprintTitle{%
 \let\@oddhead\@empty
 \let\@evenhead\@empty
 \def\@oddfoot{}%
 \let\@evenfoot\@oddfoot}
\makeatother

\usepackage{geometry}
\hypersetup{colorlinks = true, linktocpage = true, bookmarksopen = true, linkcolor = black, urlcolor=black, citecolor = black, allcolors=black}

\begin{document}

\begin{frontmatter}

{\title{{\absolutefinaledit{A fast algorithm for semi-analytically solving the homogenization boundary value problem for block locally-isotropic heterogeneous media}}}}

\author[1]{Nathan G. March\corref{cor1}}
\ead{nathan.march@hdr.qut.edu.au}

\author[1]{Elliot J. Carr}
\ead{elliot.carr@qut.edu.au}

\author[1,2]{Ian W. Turner}
\ead{i.turner@qut.edu.au}

\address[1]{School of Mathematical Sciences, Queensland University of Technology (QUT), Brisbane, Australia.}
\address[2]{ARC Centre of Excellence for Mathematical and Statistical Frontiers (ACEMS), Queensland University of Technology (QUT), Brisbane, Australia.}

\cortext[cor1]{Corresponding author}

 
\begin{abstract}
Direct numerical simulation of \newedit{diffusion} through heterogeneous media can be difficult due to the computational cost of resolving fine-scale heterogeneities. One method to overcome this difficulty is to homogenize the model by replacing the spatially-varying fine-scale diffusivity with an effective diffusivity 
\newedit{calculated} from the solution of an appropriate boundary value problem. In this paper, we present a new semi-analytical method for solving \revision{this} boundary value problem and computing the effective diffusivity for pixellated, locally-isotropic, heterogeneous media. We compare our new solution method to a standard finite volume method and show that equivalent accuracy can be achieved in less computational time for several standard test cases. We also demonstrate how the new solution method can be applied to complex heterogeneous geometries represented by \finaledit{a two-dimensional grid of rectangular blocks}. These results indicate that our new semi-analytical method has the potential to significantly speed up simulations of \newedit{diffusion} in heterogeneous \revision{media.}
\end{abstract} 


\end{frontmatter}
\section{Introduction}
\label{sec:introduction}
\finaledit{At the continuum scale, diffusion in heterogeneous media is described by the following equation:}
\begin{align}
\label{eq:heterogeneous_transport_paper2}
\frac{\partial}{\partial t}u(\mathbf{x},t) + \boldsymbol{\nabla}\cdot\left( -D(\mathbf{x})\nabla u(\mathbf{x},t)\right) = 0, \quad {\mathbf{x}} \in \Omega,
\end{align}
where \final{$D(\mathbf{x})>0$} is a \finaledit{scalar}, \amend{spatially-dependent}, isotropic diffusivity. 
\newedit{When paired with appropriate initial and boundary conditions, equation (\ref{eq:heterogeneous_transport_paper2}) is} infeasible to solve numerically if the scale at which the diffusivity $D(\mathbf{x})$ changes is small compared to the size of the domain $\Omega$, due to the prohibitively fine mesh required to capture the fine-scale heterogeneity \citep{davit_2013}. One method of overcoming this problem is to {homogenize equation (\ref{eq:heterogeneous_transport_paper2}) by replacing the \revision{spatially-varying} diffusivity $D(\mathbf{x})$ by an \final{effective} diffusivity \finaledit{tensor} $\mathbf{D}_{\mathrm{eff}}$,
yielding the homogenized equation:}
\begin{align}
\label{eq:homogeneous_transport}
\frac{\partial}{\partial t}U(\mathbf{x},t) + \boldsymbol{\nabla}\cdot\left( {-\mathbf{D}_{\mathrm{eff}}}\nabla U(\mathbf{x},t)\right) = 0, \quad {\mathbf{x}}  \in \Omega,
\end{align}
where $U(\mathbf{x},t)$ \revision{is} {a smoothed} approximation to $u(\mathbf{x},t)$. The initial and boundary conditions \newedit{for the homogenized} equation (\ref{eq:homogeneous_transport}) may take a similar form as the initial and boundary conditions \newedit{for the fine-scale} equation (\ref{eq:heterogeneous_transport_paper2}), or they may have to be modified\finaledit{, as is the case in layered geometries \citep{carr_2017a}.}
\newedit{Nevertheless, the} attraction of equation (\ref{eq:homogeneous_transport}) over equation (\ref{eq:heterogeneous_transport_paper2}) is that $\mathbf{D}_{\mathrm{eff}}$ is constant across the entire domain $\Omega$ \finalrevision{or varies across $\Omega$ at a coarse-scale compared to $D(\mathbf{x})$}. This means a coarser mesh can be used to solve the {homogenized} equation (\ref{eq:homogeneous_transport}) leading to more computationally efficient simulations. The efficiency of \newedit{solving} the {homogenized equation}, however, is negated to some extent if the overhead of computing the effective {diffusivity is high}. \finalfinaledit{For some problems, this may not be a significant issue, as the calculation of the effective diffusivities can be performed once before the simulation and stored thereafter, however the computational expense of calculating the effective diffusivities could be significant for other problems, such as those in which the geometry of $\Omega$ changes throughout the simulation.
 \finaledit{This can occur in areas such as fracturing in shale gas reservoirs \citep{zhou_2017} and cell wall collapse in wood drying \citep{carr_2013b}, as each change in the geometry \newedit{requires} the effective diffusivity to be recomputed.}} {\revision{To address this issue, the} aim of this \final{paper} is to develop an accurate and efficient method {for solving the boundary value problem required to calculate} $\mathbf{D}_{\mathrm{eff}}$}. \finaledit{While we consider two-dimensional, rectangular domains only in this work, where $\mathbf{x} = (x,y)$, the solution strategy at the core of our method can be extended to three-dimensional and circular domains}\footnote{\finalfinaledit{As we discuss in later sections of this paper, our semi-analytical method relies on representing the fluxes at the interfaces between adjacent blocks as unknown functions of space and integrating these functions. In three-dimensional, rectangular media double integrals would be required and in circular media these integrals would be computed using polar coordinates \citep{torquato_1991}}.}.

\finaledit{According to \citep{carr_2013b},} for a \finalrevision{two-dimensional} {periodic domain \revision{$\Omega$} with unit cell} \revision{$\mathcal{C} = [x_{0},x_{n}]\times [y_{0},y_{m}]$}, the first and second columns of the corresponding effective diffusivity \finaledit{tensor} can be computed using the following formulae:
\begin{gather}
\label{eq:deff_formula_1}
[\mathbf{D}_{\mathrm{eff}}]_{(:,1)} = \frac{1}{A}\int_{y_0}^{y_m}\int_{x_0}^{x_n} D(x,y) \nabla (\psi^{(x)}(x,y)+x)\, dx \, dy, \\
\label{eq:deff_formula_2}
[\mathbf{D}_{\mathrm{eff}}]_{(:,2)} = \frac{1}{A}\int_{y_0}^{y_m}\int_{x_0}^{x_n} D(x,y) \nabla (\psi^{(y)}(x,y)+y)\, dx \, dy,
\end{gather}
where \final{$x_n>x_0$, $y_m>y_0$}, $A = (x_{n}-x_{0})(y_{m}-y_{0})$\finaledit{, $\xi = x,y$} and $\psi^{(\xi)}(x,y)$ {satisfies} {the boundary value problem:}
\begin{gather}
\label{eq:diffusion}
0 = \nabla \cdot (D(x,y) \nabla(\psi^{(\xi)}\change{(x,y)}+\xi)), \quad x_0<x<x_n, \quad y_0 < y < y_m,\\
\label{eq:psi_periodic}
\text{$\psi^{(\xi)}\change{(x,y)}$ is periodic with unit cell $\mathcal{C}$,}\\
\label{eq:psi_zero_mean}
\int_{y_{0}}^{y_{m}}\int_{x_{0}}^{x_{n}} \psi^{(\xi)}(x,y)\,dx \,dy = 0.
\end{gather}
{We} note that \revision{the} formulae for the effective diffusivity (\ref{eq:deff_formula_1})--(\ref{eq:deff_formula_2}) and the boundary value problem (\ref{eq:diffusion})--(\ref{eq:psi_periodic}) are equivalent to other formulations such as those presented in \citep{dykaar_1992} and \citep{szymkiewicz_2012}.
For a layered heterogeneous medium \citep{szymkiewicz_2005}, \change{the boundary value problem} (\ref{eq:diffusion})--(\ref{eq:psi_periodic}) has an exact solution yielding an effective diffusivity consisting of area-weighted arithmetic and harmonic averages of the layer diffusivities in the directions parallel and perpendicular to the layers, respectively \citep{szymkiewicz_2005}. In general, however, the boundary value problem cannot be solved analytically \change{and one must turn to numerical schemes such as the finite volume method \citep{carr_2014} or finite element method \citep{ray_2018}}.
\begin{figure}[t]
\centering
{\includegraphics[width=\textwidth]{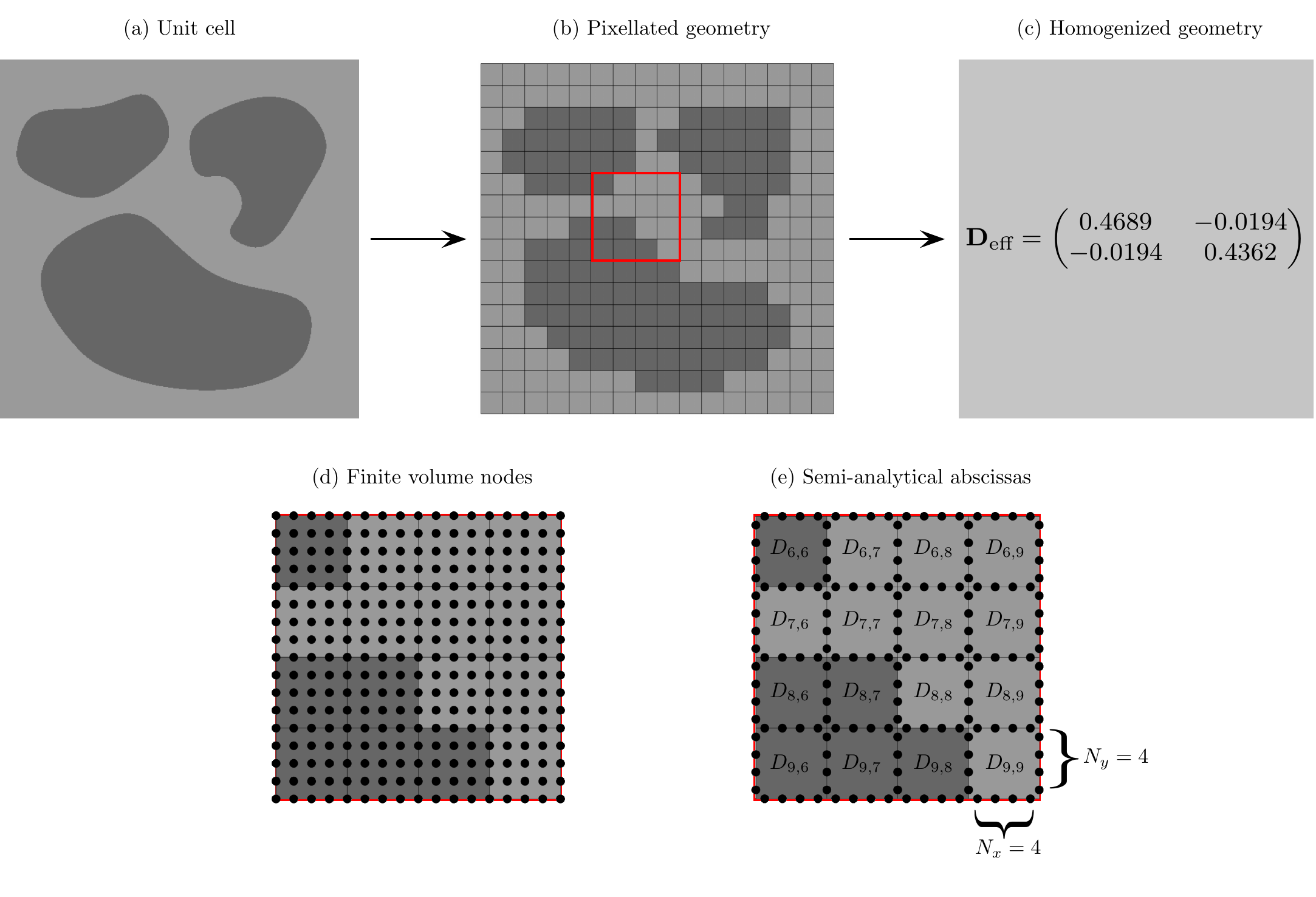}}
\caption{Homogenization of a block locally-isotropic pixelated unit cell (a) \amend{Unit cell} geometry (b) Pixelated geometry consisting of a 16 by 16 array of blocks (c) Homogenized geometry and example effective \edit{diffusivity} assuming the dark grey and light grey blocks have diffusivity 0.1 and 1.0, respectively. (d) and (e) show a zoomed in view of the section highlighted in red in (b). (d) Nodes required to solve the boundary value problem (\ref{eq:diffusion})--(\ref{eq:psi_periodic}) using the finite volume method on a structured grid (see section \ref{sec:fvm_paper2} for details) and (e) \finaledit{Diffusivities and abscissas required for our new semi-analytical method, where $D_{i,j}$ represents the diffusivity in the $(i,j)$th block and $N_x$ and $N_y$ are the number of abscissas \finalfinaledit{used to compute the integrals of functions of $x$ and $y$, respectively.}}}
\label{fig:figure1}
\end{figure}
Each of these methods \amend{require the generation of a mesh over the} \revision{unit cell $\mathcal{C}$} \edit{(see Figure \ref{fig:figure1}(d))} that resolves the fine-scale heterogeneity, \amend{which leads to a discretization of equations (\ref{eq:diffusion})--(\ref{eq:psi_periodic}) taking the form of a large linear system that is computationally expensive to solve.} \absolutefinaledit{A spectral method for \change{solving equations} (\ref{eq:diffusion})--(\ref{eq:psi_periodic}) was presented in \citep{dykaar_1992}, however they noted that these spectral methods are generally not well-suited to \change{composite media with discontinuous diffusivities}}\newedit{, such as those considered in this work,} as errors caused by Gibbs \change{phenomena} can lead to inaccurate \edit{effective diffusivities.}

{In this \final{paper}, we develop a new semi-analytical method for solving the boundary value problem (\ref{eq:diffusion})--(\ref{eq:psi_periodic})}. Our solution method is applicable \finaledit{to a heterogeneous medium taking the form of a rectangular array of blocks of any size}, where the diffusivity is \finaledit{positive,} isotropic and constant within each individual block. \finaledit{We note that these geometries are not smooth as the diffusivity can be discontinuous across the interface between adjacent blocks.} \revision{While our semi-analytical solution method is restricted to a block geometry, it can be used to homogenize complex heterogeneous geometries \newedit{by pixellating and representing the geometry as a grid of blocks, as demonstrated in Figures \ref{fig:figure1}(a)--(b).} To calculate the effective \edit{diffusivity}, our semi-analytical solution method involves reformulating \amend{the boundary value problem (\ref{eq:diffusion})--(\ref{eq:psi_periodic})} on the heterogeneous medium into a family of boundary value problems on the homogeneous individual blocks} by introducing unknown functions representing the flux at the interfaces between adjacent blocks. \absolutefinaledit{This \revision{strategy} was recently \amend{utilised} in one spatial dimension \citep{carr_2016}} and allowed accurate semi-analytical solutions to be derived for time-dependent multilayer diffusion problems involving a large number of heterogeneities (layers). In this \final{paper}, we extend these ideas to the significantly more challenging problem of steady-state diffusion \finaledit{(\ref{eq:diffusion})} in two spatial dimensions  where interfaces are no longer points but lines. \finaledit{This is achieved by solving {the boundary value problems on each homogeneous block} using standard techniques which yield analytical solution expressions that depend on integrals involving the unknown interface functions. By applying an appropriate numerical quadrature rule \edit{(see Figure \ref{fig:figure1}(e))} to these integrals and enforcing continuity of the solution across each interface, we show how the solution can be computed in each block. {The net result is a semi-analytical method for solving the well-established periodic boundary value problem (\ref{eq:diffusion})--(\ref{eq:psi_periodic}) and \revision{hence} calculating the effective diffusivity for a pixelated, locally-isotropic, heterogeneous domain}.}

%

Our new semi-analytical solution method has several desirable properties:
\begin{itemize}
\item it is able to compute, in less computational time, a solution to the boundary value problem (\ref{eq:diffusion})--(\ref{eq:psi_periodic}) to the same level of accuracy as a standard finite volume method \finaledit{(see Table \ref{tab:case1_4} for details)} \final{(the reasons for choosing a finite volume method as our benchmark method are stated in section \ref{sec:fvm_paper2})};
{\item  \finaledit{it  does not require representing the diffusivities as a function, and therefore it is not affected by the errors caused by \amend{the Gibbs phenomenon arising} in standard spectral methods \citep{dykaar_1992};}}
{\item it requires quadrature abscissas only on the interfaces and boundaries \revision{of} the \revision{unit cell} whereas standard numerical methods require a mesh across the entire \revision{unit cell} (see Figure \revision{\ref{fig:figure1}(d)--(e))}, leading to a smaller discretized linear system \revision{(see sections \ref{sec:semi-analytical_method} and \ref{sec:fvm_paper2} for details)};}
\item it can capture \revision{complex} heterogeneous {geometries, as seen in \revision{Figure \ref{fig:figure1} and later in section \ref{sec:complex_geometry}}; and
{\item it is an analytical method, which allows for the integrals \revision{defining the effective diffusivity} (\ref{eq:deff_formula_1})--(\ref{eq:deff_formula_2}) to be \revision{calculated} analytically, whereas standard numerical methods require \revision{numerical} integration to compute the effective diffusivity \revision{\citep{carr_2014}}.}}
\end{itemize}
{We stress here for clarity that our new semi-analytical \revision{method} is not a new approach to homogenize heterogeneous media but \amend{an {alternative}} \revision{for solving the boundary value problem (\ref{eq:diffusion})--(\ref{eq:psi_periodic}) and calculating the effective diffusivity (\ref{eq:deff_formula_1})--(\ref{eq:deff_formula_2}) \absolutefinaledit{\citep{hornung_1997}.} Finally, we note that our method may be considered to be a form of Trefftz method according to the definition stated in \citep{herrera_2000}, given that the domain $\Omega$ is partitioned into a number of smaller blocks and the solution across $\Omega$ is taken as a combination of solutions on these individual blocks. \finaledit{However, our method differs from Trefftz \finalrevision{methods \citep{herrera_2000}} in two key ways. Firstly, our method is fundamentally an analytical method, as we present a closed form solution to the governing PDE (\ref{eq:diffusion}) in each block and can thus evaluate the solution $\psi^{(\xi)}$ at any point in the domain $\Omega$, whereas Trefftz methods are \finalfinaledit{typically} considered to be numerical methods \citep{herrera_2000}. \neweredit{Secondly, Trefftz methods are \finalfinaledit{typically} classed as a type of finite element method \citep{herrera_2000}, which require the \finalrevision{weak} form of the governing PDE (\ref{eq:diffusion}), whereas our semi-analytical method handles the \finalrevision{strong form} directly.}}}}

\revision{The} remaining sections of this \final{paper} are organised as follows. In section \ref{sec:effective_diffusivity} we redefine the boundary value problem (\ref{eq:diffusion})--(\ref{eq:psi_periodic}) and effective diffusivity (\ref{eq:deff_formula_1})--(\ref{eq:deff_formula_2}) for a block heterogeneous medium. In section \ref{sec:semi-analytical_method}, we develop \revision{our} semi-analytical method for solving the boundary value problem \revision{(\ref{eq:diffusion})--(\ref{eq:psi_periodic})} and computing the effective diffusivity \revision{(\ref{eq:deff_formula_1})--(\ref{eq:deff_formula_2})}, describe \neweredit{the implementation of the method} and discuss the conditions required to ensure that it provides an accurate solution. In section \ref{sec:fvm_paper2}, we discuss a standard finite volume method for solving the \revision{boundary value problem and computing the effective diffusivity.} In section \ref{sec:results_paper2}, we compare the semi-analytical and finite volume methods in terms of accuracy and efficiency for a variety of test cases from the literature, apply the semi-analytical method to the homogenization of complex heterogeneous geometries, demonstrate how our semi-analytical method can be applied to complicated curved geometries by first pixelating them and \revision{discuss} \final{extensions} of our semi-analytical method \newedit{to} further increase its \revision{efficiency}. We conclude the paper in section \ref{sec:conclusion_paper2} with a summary of the key findings of the work.

\section{Definition of effective diffusivity}
\label{sec:effective_diffusivity}

Consider a rectangular heterogeneous domain \revision{$\mathcal{C} = [x_{0},x_{n}]\times [y_{0},y_{m}]$} consisting of an $m$ by $n$ array of blocks. Let the $(i,j)$th block have domain $[x_{j-1},x_{j}]\times [y_{i-1},y_{i}]$ and \final{scalar} diffusivity $D_{i,j}$ with vertical interfaces between blocks at $x = x_{j}$ ($j = 1,\hdots,n-1$) and horizontal interfaces between blocks at $y = y_{i}$ ($i = 1,\hdots,m-1$) (see Figure \ref{fig:figure1}(b)). For this domain the first and second columns of the effective diffusivity $\mathbf{D}_{\mathrm{eff}}$ are given by:
\begin{gather}
\label{eq:deff_1}
[\mathbf{D}_{\mathrm{eff}}]_{(:,1)} = \frac{1}{A}\sum_{i=1}^m\sum_{j = 1}^n\int_{y_{i-1}}^{y_i}\int_{x_{j-1}}^{x_j} D_{i,j} \nabla (\psi^{(x)}_{i,j}(x,y)+x)\, dx \, dy, \\
\label{eq:deff_1.5}
[\mathbf{D}_{\mathrm{eff}}]_{(:,2)} = \frac{1}{A}\sum_{i=1}^m\sum_{j = 1}^n\int_{y_{i-1}}^{y_i}\int_{x_{j-1}}^{x_j} D_{i,j} \nabla (\psi^{(y)}_{i,j}(x,y)+y)\, dx \, dy,
\end{gather}
where $A = (x_{n}-x_{0})(y_{m}-y_{0})$ and $\psi_{i,j}^{(\xi)}$ satisfies the following PDE on the $(i,j)$th block:
\begin{align}
\label{eq:diffusion_2}
0 = \nabla \cdot (D_{i,j}\nabla(\psi^{(\xi)}_{i,j}(x,y)+\xi)), \quad x_{j-1}<x<x_j, \quad y_{i-1} < y < y_i. 
\end{align} 
{At the interfaces between adjacent blocks both the solution and flux are continuous \citep{szymkiewicz_2005}, leading to the following interface conditions:}
\begin{gather}
\label{eq:UBC1}
\psi^{(\xi)}_{i,j}(x,y_i) =\psi^{(\xi)}_{i+1,j}(x,y_i),\\
\label{eq:UBC3}
D_{i,j}\frac{\partial}{\partial y}\left(\psi^{(\xi)}_{i,j}+\xi\right)\vline_{y = y_i} = D_{i+1,j}\frac{\partial}{\partial y}\left(\psi^{(\xi)}_{i+1,j}+\xi\right)\vline_{y = y_i},
\end{gather}
for $ x_{j-1} < x < x_j$,  $i = 1,\hdots,m-1$ and $j = 1,\hdots,n$ and
\begin{gather}
\label{eq:UBC2}
\psi^{(\xi)}_{i,j}(x_j,y) =\psi^{(\xi)}_{i,j+1}(x_j,y),\\
\label{eq:UBC4}
D_{i,j}\frac{\partial}{\partial x}\left(\psi^{(\xi)}_{i,j}+\xi\right)\vline_{x = x_j}= D_{i,j+1}\frac{\partial}{\partial x}\left(\psi^{(\xi)}_{i,j+1}+\xi\right)\vline_{x = x_j},
\end{gather}
for $ y_{i-1} < y < y_i$, $ i = 1,\hdots,m$ and $ j = 1,\hdots,n-1$,
while periodic conditions hold at the external boundaries of the domain:
\begin{gather}
\label{eq:UBC5}
\psi^{(\xi)}_{1,j}(x,y_0) =\psi^{(\xi)}_{m,j}(x,y_m),\\
\label{eq:UBC7}
D_{1,j}\frac{\partial}{\partial y}\left(\psi^{(\xi)}_{1,j}+\xi\right)\vline_{y = y_0}= D_{m,j}\frac{\partial}{\partial y}\left(\psi^{(\xi)}_{m,j}+\xi\right)\vline_{y = y_m},
\end{gather}
for $ x_{j-1} < x < x_j$ and $j = 1,\hdots,n$ and
\begin{gather}
\label{eq:UBC6}
\psi^{(\xi)}_{i,1}(x_0,y) =\psi^{(\xi)}_{i,n}(x_n,y),\\
\label{eq:UBC8}
D_{i,1}\frac{\partial}{\partial x}\left(\psi^{(\xi)}_{i,1}+\xi\right)\vline_{x = x_0}  = D_{i,n}\frac{\partial}{\partial x}\left(\psi^{(\xi)}_{i,n}+\xi\right)\vline_{x = x_n},
\end{gather}
for $y_{i-1} < y < y_i$ and $i = 1,\hdots,m.$
The boundary value problem (\ref{eq:diffusion_2})--(\ref{eq:UBC8}) has an infinite number of solutions \citep{carr_2014} as any arbitrary constant may be added to the solution $\psi_{i,j}^{(\xi)}$ and it will still satisfy the PDE (\ref{eq:diffusion_2}) as well as the boundary conditions (\ref{eq:UBC1})--(\ref{eq:UBC8}). 
\newedit{The zero mean condition (\ref{eq:psi_zero_mean}) can be rewritten in the following form:}
\begin{align}
\label{eq:zero_mean_paper2}
\frac{1}{A}\sum_{i=1}^m\sum_{j = 1}^n \int_{y_{i-1}}^{y_i}\int_{x_{j-1}}^{x_j} \psi^{(\xi)}_{i,j}(x,y)\,dx\,dy = 0.
\end{align}
With $\psi_{i,j}^{(\xi)}(x,y)$ known for $\xi = x$ and $\xi = y$, the effective diffusivity, $\mathbf{D}_{\mathrm{eff}}$, can be computed from equations (\ref{eq:deff_1})--(\ref{eq:deff_1.5}). Interpretation of the entries of the symmetric positive definite matrix $\mathbf{D}_{\mathrm{eff}}$ comes from its diagonalization, $\mathbf{D}_{\mathrm{eff}}= \mathbf{P}\boldsymbol{\Lambda}\mathbf{P}^T$. Here, $\mathbf{P} \in \mathbb{R}^{2 \times 2}$ is an orthonormal matrix whose columns are unit vectors pointing in the principal directions of diffusion and $\boldsymbol{\Lambda} \in \mathbb{R}^{2 \times 2}$ is a diagonal matrix whose entries represent the apparent diffusivities in these directions \citep{plawsky_2009}.

\section{Semi-analytical method}
\label{sec:semi-analytical_method}
\subsection{Solution approach}
\begin{figure}[t]
\centering
{\includegraphics[width=\textwidth]{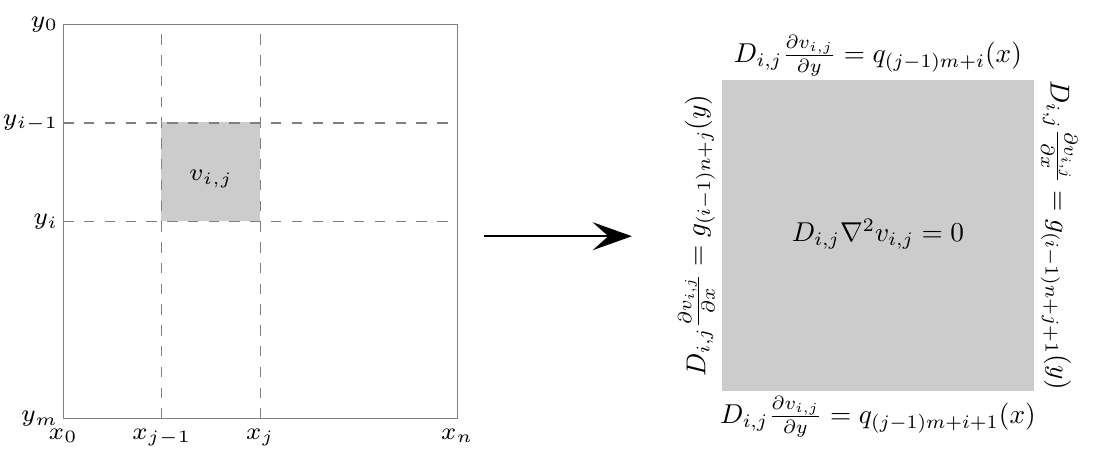}}
\caption{A subsection of the domain displaying the $(i,j)$th block (left) and the boundary value problem (\ref{eq:diffusion_4})--(\ref{eq:BC2}) imposed on the $(i,j)$th block (right).}
\label{fig:individual_block}
\end{figure}

We now describe our new semi-analytical solution approach for solving the boundary value problem (\ref{eq:diffusion_2})--(\ref{eq:zero_mean_paper2}) which enables the effective diffusivity (\ref{eq:deff_1})-(\ref{eq:deff_1.5}) to be determined. We first transform equation (\ref{eq:diffusion_2}) into Laplace's equation via the transformation:
\begin{align}
\label{eq:transformation}
v_{i,j}^{(\xi)}(x,y) = \psi^{(\xi)}_{i,j}(x,y) + \xi,
\end{align}
for $i = 1,\hdots,m$ and $j = 1,\hdots,n$. \finaledit{Combining equation (\ref{eq:transformation}) with equations (\ref{eq:diffusion_2}) yields:}
\begin{align}
\label{eq:diffusion_3}
D_{i,j}\left(\frac{\partial^2 v_{i,j}^{(\xi)}}{\partial x^2} + \frac{\partial^2 v_{i,j}^{(\xi)}}{\partial y^2}\right) = 0, 
\end{align}
for $ x_{j-1}<x<x_j$ and $y_{i-1} < y < y_i$ and \finaledit{combining equation (\ref{eq:transformation}) with equations (\ref{eq:UBC1}), (\ref{eq:UBC2}), (\ref{eq:UBC5}) and (\ref{eq:UBC6}) yields the following conditions on the solution $v_{i,j}^{(\xi)}$ at the interfaces and boundaries:}
\begin{align}
\label{eq:VBC1}
v_{i,j}^{(\xi)}(x,y_i) &=v_{i+1,j}^{(\xi)}(x,y_i), \quad i = 1,\hdots,m-1, \quad j = 1,\hdots,n,\\
\label{eq:VBC2}
v_{i,j}^{(\xi)}(x_j,y) &=v_{i,j+1}^{(\xi)}(x_j,y), \quad i = 1,\hdots,m, \quad j = 1,\hdots,n-1,\\
\label{eq:VBC3}
v_{i,1}^{(\xi)}(x_0,y) &= \begin{cases} v_{i,n}^{(\xi)}(x_n,y)+x_0-x_n, \quad &\text{if} \quad \xi = x, \\ 
v_{i,n}^{(\xi)}(x_n,y), \quad &\text{if} \quad \xi = y, \end{cases} \quad i = 1,\hdots,m,\\
\label{eq:VBC4}
v_{1,j}^{(\xi)}(x,y_0) &= \begin{cases} v_{m,j}^{(\xi)}(x,y_m), \quad &\text{if} \quad \xi = x, \\ 
v_{m,j}^{(\xi)}(x,y_m)+y_0-y_m, \quad &\text{if} \quad \xi = y, \end{cases} \quad j = 1,\hdots,n,
\end{align}
for $ x_{j-1}<x<x_j$ and $y_{i-1} < y < y_i$ and \finaledit{combining equation (\ref{eq:transformation}) with equations (\ref{eq:UBC3}), (\ref{eq:UBC4}), (\ref{eq:UBC7}) and (\ref{eq:UBC8}) yields the following conditions on the flux at the interfaces and boundaries:}
\begin{gather}
\label{eq:VBC9}
D_{i,1}\frac{\partial v_{i,1}^{(\xi)}}{\partial x}(x_0,y) =D_{i,n}\frac{\partial v_{i,n}^{(\xi)}}{\partial y}(x_n,y), \quad i = 1,\hdots,m,\\
\label{eq:VBC10}
D_{1,j}\frac{\partial v_{1,j}^{(\xi)}}{\partial y}(x,y_0) =D_{m,j}\frac{\partial v_{m,j}^{(\xi)}}{\partial y}(x,y_m), \quad j = 1,\hdots,n,\\
\label{eq:VBC11}
D_{i,j}\frac{\partial v_{i,j}^{(\xi)}}{\partial y}(x,y_i) =D_{i+1,j}\frac{\partial v_{i+1,j}^{(\xi)}}{\partial y}(x,y_i),  \quad i = 1,\hdots,m-1, \quad j = 1,\hdots,n,\\
\label{eq:VBC12}
D_{i,j}\frac{\partial v_{i,j}^{(\xi)}}{\partial x}(x_j,y) =D_{i,j+1}\frac{\partial v_{i,j+1}^{(\xi)}}{\partial x}(x_j,y), \quad i = 1,\hdots,m, \quad j = 1,\hdots,n-1,
\end{gather}
for $ x_{j-1}<x<x_j$ and $y_{i-1} < y < y_i$. \finaledit{Combining the zero mean condition (\ref{eq:zero_mean_paper2}) with the transformation (\ref{eq:transformation}) yields:}
\begin{align}
\label{eq:zero_mean_2}
\sum_{i=1}^m\sum_{j = 1}^n \int_{y_{i-1}}^{y_i}\int_{x_{j-1}}^{x_j} v_{i,j}^{(\xi)}(x,y)\,dx\,dy &= \begin{cases} \displaystyle{\frac{x_n+x_0}{2}}, \quad \text{if} \quad \xi = x,\vspace*{0.1cm} \\ 
\displaystyle{\frac{y_m+y_0}{2}}, \quad \text{if} \quad \xi = y. \end{cases}
\end{align}
Our semi-analytical approach for solving this transformed boundary value problem involves first setting the fluxes at the interfaces and boundaries in equations (\ref{eq:VBC9})--(\ref{eq:VBC12}) to be unknown \final{spatial functions} \citep{carr_2018a} as follows:
\begin{gather}
\label{eq:flux1}
D_{i,1}\frac{\partial v_{i,1}^{(\xi)}}{\partial x}(x_0,y) =D_{i,n}\frac{\partial v_{i,n}^{(\xi)}}{\partial y}(x_n,y):=g_{(i-1)n+1}(y), \quad i = 1,\hdots,m,\\
\label{eq:flux2}
D_{1,j}\frac{\partial v_{1,j}^{(\xi)}}{\partial y}(x,y_0) =D_{m,j}\frac{\partial v_{m,j}^{(\xi)}}{\partial y}(x,y_m):=q_{(j-1)m+1}(x), \quad j = 1,\hdots,n,\\
\label{eq:flux3}
\finalrevision{D_{i,j}\frac{\partial v_{i,j}^{(\xi)}}{\partial x}(x_j,y) =D_{i,j+1}\frac{\partial v_{i,j+1}^{(\xi)}}{\partial x}(x_j,y):=
g_{(i-1)n+j+1}(y)}, \quad i = 1,\hdots,m, \quad j = 1,\hdots,n-1,\\
\label{eq:flux4}
\finalrevision{D_{i,j}\frac{\partial v_{i,j}^{(\xi)}}{\partial y}(x,y_i) =D_{i+1,j}\frac{\partial v_{i+1,j}^{(\xi)}}{\partial y}(x,y_i):=
q_{(j-1)m+i+1}(x)}, \quad i = 1,\hdots,m-1, \quad j = 1,\hdots,n\finalrevision{,}
\end{gather}
for $ x_{j-1}<x<x_j$ and $y_{i-1} < y < y_i$. This allows us to reformulate the heterogeneous boundary value problem (\ref{eq:diffusion_3})--(\ref{eq:VBC12}) into a family of boundary value problems on each of the $mn$ blocks \finaledit{(see Figure \ref{fig:individual_block}):}
\begin{gather}
\label{eq:diffusion_4}
D_{i,j}\left(\frac{\partial^2 v_{i,j}^{(\xi)}}{\partial x^2} + \frac{\partial^2 v_{i,j}^{(\xi)}}{\partial y^2}\right) = 0, \quad x_{j-1}<x<x_j, \quad y_{i-1} < y < y_i,\\
\label{eq:BC1}
D_{i,j}\frac{\partial v_{i,j}^{(\xi)}}{\partial x}(x_{j-1},y) = g_{(i-1)n+j}(y),\quad
D_{i,j}\frac{\partial v_{i,j}^{(\xi)}}{\partial x}(x_j,y) = g_{(i-1)n+j+1}(y), \quad y_{i-1} < y < y_i,\\
D_{i,j}\frac{\partial v_{i,j}^{(\xi)}}{\partial y}(x,y_{i-1}) =q_{(j-1)m+i}(x),\quad
\label{eq:BC2}
D_{i,j}\frac{\partial v_{i,j}^{(\xi)}}{\partial y}(x,y_i) =q_{(j-1)m+i+1}(x), \quad x_{j-1} < x < x_j,
\end{gather}
where the unknown interface and boundary functions, introduced in equations (\ref{eq:flux1})--(\ref{eq:flux4}), are constrained to satisfy the interface and boundary conditions (\ref{eq:VBC1})--(\ref{eq:VBC4}). The boundary conditions (\ref{eq:BC1})--(\ref{eq:BC2}) are valid for all blocks except those in either the bottom row ($i = m$) {or} right column $(j = n)$. The boundary conditions for these blocks are similar to equations (\ref{eq:BC1})--(\ref{eq:BC2}) and detailed in \ref{sec:BCs}.

The solution of the boundary value problem (\ref{eq:diffusion_4})--(\ref{eq:BC2}) is given by \citep{polyanin_2002}\footnote{\finalrevision{For three-dimensional problems, the equivalent solution would be more complicated, requiring a double summation and the coefficients to be computed in terms of double integrals.}}:
\begin{equation}
\begin{aligned}
\label{eq:solution}
v_{i,j}^{(\xi)}(x,y) &= -\frac{a_{i,j,0}^{(\xi)}}{4l_j}(x-x_j)^2 + \frac{b_{i,j,0}^{(\xi)}}{4l_j}(x-x_{j-1})^2 - \frac{c_{i,j,0}^{(\xi)}}{4h_i}(y-y_i)^2 + \frac{d_{i,j,0}^{(\xi)}}{4h_i}(y-y_{i-1})^2+K_{i,j}\\
&-h_i \sum_{k = 1}^\infty  \frac{a_{i,j,k}^{(\xi)}}{\gamma_{i,j,k}} \cosh\left[\frac{k\pi(x-x_j)}{h_i}\right]\cos\left(\frac{k\pi(y-y_{i-1})}{h_i}\right)\\
&+h_i \sum_{k = 1}^\infty  \frac{b_{i,j,k}^{(\xi)}}{\gamma_{i,j,k}} \cosh\left[\frac{k\pi(x-x_{j-1})}{h_i}\right]\cos\left(\frac{k\pi(y-y_{i-1})}{h_i}\right) \\ 
&-l_j\sum_{k = 1}^\infty  \frac{c_{i,j,k}^{(\xi)}}{\mu_{i,j,k}} \cosh\left[\frac{k\pi(y-y_i)}{l_j}\right]\cos\left(\frac{k\pi(x-x_{j-1})}{l_j}\right)\\
&+ l_j\sum_{k = 1}^\infty  \frac{d_{i,j,k}^{(\xi)}}{\mu_{i,j,k}} \cosh\left[\frac{k\pi(y-y_{i-1})}{l_j}\right]\cos\left(\frac{k\pi(x-x_{j-1})}{l_j}\right),
\end{aligned}
\end{equation}
where:
\begin{align}
h_i &= y_i - y_{i-1}, \quad l_j = x_j - x_{j-1},\\
\gamma_{i,j,k}&= k\pi\sinh\frac{k\pi l_j}{h_i},\quad
\mu_{i,j,k}= k\pi\sinh\frac{k\pi h_i}{l_j},
\end{align}
and $K_{i,j}$ is an arbitrary constant. While for the solution of the boundary value problem (\ref{eq:diffusion_4})--(\ref{eq:BC2}) the constant $K_{i,j}$ is arbitrary, in the context of the solution across the entire \revision{unit cell $\mathcal{C}$} consisting of all $mn$ boundary value problems, the constants are not arbitrary as they are required to ensure continuity of the solution across the interfaces and boundaries. Note that the coefficients $a_{i,j,k}^{(\xi)}$, $b_{i,j,k}^{(\xi)}$, $c_{i,j,k}^{(\xi)}$ and $d_{i,j,k}^{(\xi)}$ appearing in equation (\ref{eq:solution}) depend on the values of $i$ and $j$. For $i = 1,\hdots,m-1$ and $j = 1,\hdots,n-1$:
\begin{align}
\label{eq:coeff_1_paper2}
a_{i,j,k}^{(\xi)}&= \frac{2}{h_i} \int_{y_{i-1}}^{y_i} \frac{g_{(i-1)n+j}(y)}{D_{i,j}} \cos\left(\frac{k\pi (y-y_{i-1})}{h_i}\right)\, dy,\\
\label{eq:coeff_2_paper2}
b_{i,j,k}^{(\xi)}&= \frac{2}{h_i} \int_{y_{i-1}}^{y_i} \frac{g_{(i-1)n+j+1}(y)}{D_{i,j}} \cos\left(\frac{k\pi (y-y_{i-1})}{h_i}\right)\, dy,\\
\label{eq:coeff_3_paper2}
c_{i,j,k}^{(\xi)}&= \frac{2}{l_j} \int_{x_{j-1}}^{x_j} \frac{q_{(j-1)m+i}(x)}{D_{i,j}} \cos\left(\frac{k\pi(x-x_{j-1})}{l_j}\right)\, dx,\\
\label{eq:coeff_4_paper2}
d_{i,j,k}^{(\xi)}&= \frac{2}{l_j} \int_{x_{j-1}}^{x_j} \frac{q_{(j-1)m+i+1}(x)}{D_{i,j}} \cos\left(\frac{k\pi (x-x_{j-1})}{l_j}\right)\, dx,
\end{align}
{where $k = 0,1,\hdots,\infty$,} while for the bottom row ($i = m$) and right column ($j = n$) of blocks the coefficients are similar and given in \ref{sec:coefficients}.

Due to the boundary conditions (\ref{eq:BC1})--(\ref{eq:BC2}), the solution (\ref{eq:solution}) requires an additional solvability condition to ensure a net flux of zero is imposed on each block \citep{polyanin_2002}. This condition takes the form:
\begin{multline}
\label{eq:S1}
\int_{y_{i-1}}^{y_i} \frac{g_{(i-1)n+j}(y)}{D_{i,j}}\, dy - \int_{y_{i-1}}^{y_i} \frac{g_{(i-1)n+j+1}(y)}{D_{i,j}}\, dy \\+  \int_{x_{j-1}}^{x_j} \frac{q_{(j-1)m+i}(x)}{D_{i,j}} \, dx - \int_{x_{j-1}}^{x_j} \frac{q_{(j-1)m+i+1}(x)}{D_{i,j}} \, dx = 0,
\end{multline}
for $i = 1,\hdots,m-1$ and $j = 1,\hdots,n-1$. The solvability conditions for the bottom row ($i = m$) and right column ($j = n$) of blocks have a similar form to equation (\ref{eq:S1}) and are given in \finalrevision{\ref{sec:solvability}.} 

Through the coefficients (\ref{eq:coeff_1_paper2})--(\ref{eq:coeff_4_paper2}), the solution (\ref{eq:solution}) in each block involves the integral of the unknown interface and boundary functions, introduced in equations  (\ref{eq:flux1})--(\ref{eq:flux4}), and the unknown constants $K_{i,j}$. By applying the as yet unused interface and boundary conditions (\ref{eq:VBC1})--(\ref{eq:VBC4}), we now derive a linear system of equations {whose solution} allows the integrals appearing in the coefficients (\ref{eq:coeff_1_paper2})--(\ref{eq:coeff_4_paper2}) to be approximated. {This is achieved by first applying a quadrature rule to these integrals}. For example using the expressions for $a_{i,j,k}^{(\xi)}$ and $c_{i,j,k}^{(\xi)}$ in equations (\ref{eq:coeff_1_paper2})--(\ref{eq:coeff_4_paper2}), we have:
\begin{multline}
\label{eq:quadrature}
a_{i,j,k}^{(\xi)}= \frac{2}{h_i} \int_{y_{i-1}}^{y_i} \frac{g_{(i-1)n+j}(y)}{D_{i,j}} \cos\left(\frac{k\pi (y-y_{i-1})}{h_i}\right)\, dy \\\approx \frac{2}{D_{i,j}h_i}\sum_{p = 1}^{N_y} w_y^{(i,p)} g_{(i-1)n+j}(y^{(i,p)}) \cos\left(\frac{k\pi (y^{(i,p)}-y_{i-1})}{h_i}\right),
\end{multline}
\begin{multline}
\label{eq:quadrature2}
\hspace*{-0.1365cm}c_{i,j,k}^{(\xi)}= \frac{2}{l_j} \int_{x_{j-1}}^{x_j} \frac{q_{(j-1)m+i}(x)}{D_{i,j}} \cos\left(\frac{k\pi(x-x_{j-1})}{l_j}\right)\, dx\\ \approx \frac{2}{D_{i,j}l_j}\sum_{p = 1}^{N_x} w_x^{(j,p)} q_{(j-1)m+i}(x^{(j,p)}) \cos\left(\frac{k\pi (x^{(j,p)}-x_{j-1})}{l_j}\right),
\end{multline}
\finalrevision{where} $w_x^{(j,p)}$ and $w_y^{(i,p)}$ are the weights of the $p$th {abscissas, $x^{(j,p)}$ and $y^{(i,p)}$,} respectively and $N_x$ and $N_y$ denote the number of abscissas. The weights $w_x^{(j,p)}$ and $w_y^{(i,p)}$ and abscissas $x^{(j,p)}$ and $y^{(i,p)}$ depend on the type of quadrature rule chosen, as will be discussed in the following section. {There are a total of $2mn$ unknown flux functions, $q_{(j-1)m+i}(x)$ and $g_{(i-1)n+j}(y)$ for $i = 1,\hdots,m$ and $j = 1,\hdots,n$, with the integrals involving these functions approximated using a quadrature rule with $N_x$ and $N_y$ function evaluations, respectively. In total, this yields $mn(N_x+N_y)$ unknown evaluations at the abscissas: $q_{(j-1)m+i}(x^{(j,p)})$ for $i = 1,\hdots,m$, $j = 1,\hdots,n$ and $p = 1,\hdots,N_{x}$ and $g_{(i-1)n+j}(y^{(i,p)})$ for $i = 1,\hdots,m$, $j = 1,\hdots,n$ and $p = 1,\hdots,N_{y}$.} To compute these evaluations, which appear in equations (\ref{eq:quadrature})--(\ref{eq:quadrature2}), we make use of the boundary and interface conditions (\ref{eq:VBC1})--(\ref{eq:VBC4}) evaluated at each of the abscissas: 
\begin{align}
\label{eq:ICX1}
v_{i+1,j}^{(\xi)}(x^{(j,p)},y_i)-v_{i,j}^{(\xi)}(x^{(j,p)},y_i) &= 0,  \quad i = 1,\hdots,m-1, \quad j = 1,\hdots,n, \\
\label{eq:ICX3}
v_{m,j}^{(\xi)}(x^{(j,p)},y_m) - v_{1,j}^{(\xi)}(x^{(j,p)},y_0) &= \begin{cases} 0, \quad &\text{if} \quad \xi = x, \\  y_m-y_0, \quad &\text{if} \quad \xi = y, \end{cases} \quad j = 1,\hdots,n,
\end{align}
for $x_{j-1} < x^{(j,p)}  <  x_j$ and $p = 1,\hdots,N_x$, 
\begin{align}
\label{eq:ICX2}
v_{i,j+1}^{(\xi)}(x_j,y^{(i,p)}) - v_{i,j}^{(\xi)}(x_j,y^{(i,p)}) &= 0, \quad i = 1,\hdots,m, \quad j = 1,\hdots,n-1,\\
\label{eq:ICX4}
v_{i,n}^{(\xi)}(x_n,y^{(i,p)})-v_{i,1}^{(\xi)}(x_0,y^{(i,p)}) &= \begin{cases} x_n-x_0, \quad &\text{if} \quad \xi = x, \\  0, \quad &\text{if} \quad \xi = y, \end{cases} \quad i = 1,\hdots,m,
\end{align}
for $y_{i-1}  <  y^{(i,p)}  <  y_i$ and $p = 1,\hdots,N_y$.
{Equations (\ref{eq:ICX1})--(\ref{eq:ICX4}) define a system of $mn(N_x+N_y)$ equations linear in the $mn(N_x+N_y)$ unknown evaluations at the abscissas. In addition to these unknowns, we have an additional $mn$ unknowns corresponding to the constants $K_{i,j}$ appearing in equation (\ref{eq:solution}). These additional unknowns are accounted for by imposing} the solvability conditions (\ref{eq:S1}), with the integrals replaced with quadrature approximations similar to {those used} in equations (\ref{eq:quadrature})--(\ref{eq:quadrature2}). The solvability conditions ensure a zero net flux on each block and as periodic conditions are imposed on the external boundaries of the \revision{unit cell $\mathcal{C}$}, the solvability condition in any one of the $mn$ blocks is redundant. To {account for this redundancy}, we replace the solvability conditions on block $mn$ with the zero mean condition (\ref{eq:zero_mean_paper2}). In the solution expression (\ref{eq:solution}), all terms with the exception of the constants $K_{i,j}$ integrate to zero over their domain of integration $[x_{j-1},x_j] \times [y_{i-1},y_i]$. This allows for the zero mean condition (\ref{eq:zero_mean_2}) to be written in the following form:
\begin{align}
\label{eq:zero_mean_K}
\sum_{i=1}^m\sum_{j = 1}^n K_{i,j} 
 &= \begin{cases} \displaystyle{\frac{x_n+x_0}{2}}, \quad \text{if} \quad \xi = x, \vspace*{0.1cm}\\ 
\displaystyle{\frac{y_m+y_0}{2}}, \quad \text{if} \quad \xi = y. \end{cases}
\end{align}
Assembling equations (\ref{eq:ICX1})--(\ref{eq:zero_mean_K}) yields a linear system:
\begin{align}
\label{eq:linear_system_paper2}
\mathbf{A}_{\text{S}}\mathbf{x}_{\text{S}} = \mathbf{b}_{\text{S}},
\end{align}
where $\mathbf{A}_{\text{S}} \in \mathbb{R}^{N \times N}$,  {$\mathbf{x}_{\text{S}} \in \mathbb{R}^{N}$, $\mathbf{b}_{\text{S}} \in \mathbb{R}^{N}$} and $N = mn(N_x+N_y+1)$.  The entries of $\mathbf{A}_{\text{S}}$ and $\mathbf{b}_{\text{S}}$ are determined via the solvability conditions (\ref{eq:S1}), the continuity conditions (\ref{eq:ICX1})--(\ref{eq:ICX4}) and {the} zero mean condition (\ref{eq:zero_mean_K}). The vector $\mathbf{x}_{\text{S}}$ contains the function evaluations {$q_{(j-1)m+i}(x^{(j,p)})$ for $i = 1,\hdots,m$, $j = 1,\hdots,n$ and $p = 1,\hdots,N_{x}$ and $g_{(i-1)n+j}(y^{(i,p)})$ for $i = 1,\hdots,m$, $j = 1,\hdots,n$ and $p = 1,\hdots,N_{y}$ as well as the constants $K_{i,j}$ for $i = 1,\hdots,m$ and $j = 1,\hdots,n$}. {Note that when constructing the linear system, we truncate the summations appearing in the solution expression (\ref{eq:solution}) after {a finite number of terms}, $k = N_\text{eig}$}. Solving the linear system  (\ref{eq:linear_system_paper2}) then allows the approximations to the integrals (\ref{eq:quadrature}) and (\ref{eq:quadrature2}) to be computed and thus the solution (\ref{eq:solution}) can be evaluated. With $v_{i,j}^{(\xi)}(x,y)$ computed, the solution of the original boundary value problem (\ref{eq:diffusion_2})--(\ref{eq:zero_mean_paper2}) can be calculated as $\psi^{(\xi)}_{i,j}(x,y) = v_{i,j}^{(\xi)}(x,y) - \xi$. 
 
In terms of the transformed solution $v_{i,j}^{(\xi)}(x,y)$, the four entries of the {effective diffusivity} are given by: 
\begin{gather}
\label{eq:diff_1}
[\mathbf{D}_{\mathrm{eff}}]_{(1,1)} = \frac{1}{A} \sum_{i = 1}^m \sum_{j = 1}^n  \int_{y_{i-1}}^{y_i}\int_{x_{j-1}}^{x_j} D_{i,j} \frac{\partial v^{(x)}_{i,j}}{\partial x}\, dx \, dy, \\
[\mathbf{D}_{\mathrm{eff}}]_{(2,1)} =\frac{1}{A}  \sum_{i = 1}^m \sum_{j = 1}^n \int_{y_{i-1}}^{y_i}\int_{x_{j-1}}^{x_j} D_{i,j} \frac{\partial v^{(x)}_{i,j}}{\partial y}\, dx \, dy,\\
[\mathbf{D}_{\mathrm{eff}}]_{(1,2)} = \frac{1}{A}  \sum_{i = 1}^m \sum_{j = 1}^n \int_{y_{i-1}}^{y_i}\int_{x_{j-1}}^{x_j} D_{i,j} \frac{\partial v^{(y)}_{i,j}}{\partial x}\, dx \, dy,\\
[\mathbf{D}_{\mathrm{eff}}]_{(2,2)} = \frac{1}{A} \sum_{i = 1}^m \sum_{j = 1}^n \int_{y_{i-1}}^{y_i}\int_{x_{j-1}}^{x_j} D_{i,j} \frac{\partial v^{(y)}_{i,j}}{\partial y}\, dx \, dy.
\label{eq:diff_4}
\end{gather}
As we have an {analytical expression for $v_{i,j}^{(\xi)}(x,y)$ (\ref{eq:solution})}{,} we can {evaluate the derivatives and integrals} appearing in equations (\ref{eq:diff_1})--(\ref{eq:diff_4}) {directly}, yielding:
\begin{gather}
\label{eq:eff_diff_1}
\hspace{-0.3cm}[\mathbf{D}_{\mathrm{eff}}]_{(1,1)} = \frac{1}{A}  \sum_{i = 1}^m \sum_{j = 1}^n \frac{D_{i,j}A_{i,j} (a_{i,j,0}^{(x)} + b_{i,j,0}^{(x)})}{4}+l_j^2 \sum_{k = 1}^{N_\text{eig}}  \frac{(c_{i,j,k}^{(x)}-d_{i,j,k}^{(x)})[1-(-1)^k]}{k\pi},\\
\hspace{-0.3cm}[\mathbf{D}_{\mathrm{eff}}]_{(1,2)} = \frac{1}{A}  \sum_{i = 1}^m \sum_{j = 1}^n \frac{D_{i,j}A_{i,j} (a_{i,j,0}^{(y)} + b_{i,j,0}^{(y)})}{4}+l_j^2 \sum_{k = 1}^{N_\text{eig}}  \frac{(c_{i,j,k}^{(y)}-d_{i,j,k}^{(y)})[1-(-1)^k]}{k\pi},\\
\hspace{-0.3cm}[\mathbf{D}_{\mathrm{eff}}]_{(2,1)} = \frac{1}{A} \sum_{i = 1}^m \sum_{j = 1}^n  \frac{D_{i,j}A_{i,j}(c_{i,j,0}^{(x)}+d_{i,j,0}^{(x)})}{4}+ h_i^2\sum_{k = 1}^{N_\text{eig}}  \frac{(a_{i,j,k}^{(x)}-b_{i,j,k}^{(x)})[1-(-1)^k]}{k\pi},\\
\label{eq:eff_diff_4}
\hspace{-0.3cm}[\mathbf{D}_{\mathrm{eff}}]_{(2,2)} = \frac{1}{A}  \sum_{i = 1}^m \sum_{j = 1}^n  \frac{D_{i,j}A_{i,j}(c_{i,j,0}^{(y)}+d_{i,j,0}^{(y)})}{4}+ h_i^2\sum_{k = 1}^{N_\text{eig}}  \frac{(a_{i,j,k}^{(y)}-b_{i,j,k}^{(y)})[1-(-1)^k]}{k\pi}\finaledit{,}
\end{gather}
\finaledit{where $A_{i,j} = (x_j - x_{j-1})(y_i - y_{i-1})$.} Finally, we remark that the coefficient matrix $\mathbf{A}_{\text{S}}$ in {the linear system} (\ref{eq:linear_system_paper2}) is identical for both $\xi = x$ and $\xi = y$ \final{, however the vector $\mathbf{b}_{\text{S}}$ differs}. This means that the {linear system for $\xi = x$ and the linear system for $\xi = y$} can both be solved with only one matrix factorisation.
\subsection{Choice of quadrature method}
\label{sec:implementation}
There are several constraints governing {the} choice of quadrature method utilized in the approximations of the integrals (\ref{eq:quadrature})--(\ref{eq:quadrature2}){:}
\begin{itemize}
\item \textit{Abscissas cannot include the vertices of blocks.}\\
If a vertex of a block is used as an abscissa in the quadrature rules (\ref{eq:quadrature})--(\ref{eq:quadrature2}), then the set of four continuity conditions (\ref{eq:ICX1})--(\ref{eq:ICX4}) evaluated at the vertex will be linearly dependent. If all $mn$ vertices are used in the quadrature rules (\ref{eq:quadrature})--(\ref{eq:quadrature2}), then the linear system (\ref{eq:linear_system_paper2}) is deficient in rank by $mn$. As the vertices of the blocks are the limits of integration for all integrals appearing in the coefficients (\ref{eq:coeff_1_paper2})--(\ref{eq:coeff_4_paper2}), to ensure linear independence we {implement} a quadrature rule that uses at most one of the limits of integration as an abscissa. While we could implement different quadrature rules that use endpoints for some integrals and do not use endpoints for other integrals, we take the \neweredit{consistent} option of using no endpoints.
\item \textit{The quadrature rule must use the same abscissas for all frequencies $k = 0,\hdots,N_{\text{eig}}$ appearing in the coefficients (\ref{eq:coeff_1_paper2})--(\ref{eq:coeff_4_paper2}).}\\
The linear system (\ref{eq:linear_system_paper2}) is generated by evaluating the boundary and interface conditions (\ref{eq:ICX1})--(\ref{eq:ICX4}) at the abscissas used in the quadrature formulae (\ref{eq:quadrature})-(\ref{eq:quadrature2}). In order to minimise the size of the the linear system (\ref{eq:linear_system_paper2}), the same abscissas $x^{(j,p)}$ and $y^{(i,p)}$ should be used for all $k = 0,\hdots,N_{\text{eig}}$. Hence, we restrict the choice of quadrature rules \neweredit{to those that use the same abscissas for all frequencies $k = 0,\hdots,N_{\text{eig}}$.}
\neweredit{\item  \textit{The quadrature rule must be able to compute the integrals appearing in equations (\ref{eq:coeff_1_paper2})--(\ref{eq:coeff_4_paper2}) for all frequencies $k = 0,\hdots,N_{\text{eig}}$ for $N_{\text{eig}} = \text{max}(N_x,N_y) - 1$ to a sufficient level of accuracy.}\\
The boundary and interface conditions (\ref{eq:ICX1})--(\ref{eq:ICX4}) result in the coefficient matrix $\mathbf{A}_{S}$ in the linear system (\ref{eq:linear_system_paper2}) being rank deficient if $N_{\text{eig}} = \text{max}(N_x,N_y) - 1$. As there is a lower limit on $N_{\text{eig}}$, the quadrature rule must be able to accurately compute the integrals appearing in equations (\ref{eq:coeff_1_paper2})--(\ref{eq:coeff_4_paper2}) for all frequencies $k = 0,\hdots,N_{\text{eig}}$, for $N_{\text{eig}} = \text{max}(N_x,N_y) - 1$.} 

\finaledit{\item \textit{The quadrature rule should be more accurate for low frequencies than for high frequencies.}\\
The entries in the coefficient matrix $\mathbf{A}_{S}$ appearing in the linear system (\ref{eq:linear_system_paper2}) are generated through a summation from $k = 0,\hdots,N_{\text{eig}}$ of evaluations of the solution (\ref{eq:solution}) at the abscissas used in the quadrature rule. As $k$ increases, these solution evaluations decrease and therefore solution evaluations with low frequencies comprise a larger proportion of the \finalfinaledit{entries of the }coefficient matrix $\mathbf{A}_{S}$ than solution evaluations with higher frequencies. Therefore, the quadrature rule should be more accurate for low frequencies than for high frequencies.}

\item \neweredit{\textit{The quadrature rule must be a linear function of the function evaluations and cannot make use of any information about the function being integrated}.}\\
As the evaluations, {$q_{(j-1)m+i}(x^{(j,p)})$ for $i = 1,\hdots,m$, $j = 1,\hdots,n$ and $p = 1,\hdots,N_{x}$ and $g_{(i-1)n+j}(y^{(i,p)})$ for $i = 1,\hdots,m$, $j = 1,\hdots,n$ and $p = 1,\hdots,N_{y}$, used to approximate the integrals (\ref{eq:coeff_1_paper2})--(\ref{eq:coeff_4_paper2}) are} determined via the solution of the linear system (\ref{eq:linear_system_paper2}), the quadrature rule must be a linear function of these evaluations and cannot use any other information about the function, \neweredit{such as its derivatives, which are commonly used in approximations of oscillatory integrals \citep{dominguez_2014}}.
\end{itemize}

\finaledit{As the integrals appearing in the coefficients (\ref{eq:coeff_1_paper2})--(\ref{eq:coeff_4_paper2}) contain an oscillatory function, we investigated Filon \citep{filon_1928} and Levin \citep{levin_1996} type methods, which are quadrature rules well-suited to highly oscillatory functions. The choice of quadrature rule is further complicated by the integrands being unknown. We investigated a wide variety of different Filon and Levin type quadrature methods, including (but not limited to) the many different methods presented in \citep{deano_2018} and found that almost all of them were not suitable for use in our semi-analytical method, as they do not satisfy one or more of the above mentioned constraints. The quadrature methods that were found to be  suitable are: the midpoint rule, a Filon-midpoint rule demonstrated in \citep{potticary_2005} and various Gaussian quadrature methods that do not use abscissas at the endpoints. } {The Filon-midpoint rule was found to be consistently less accurate than the regular midpoint rule in \neweredit{testing its application to oscillatory integrals}. We then tested our semi-analytical method with both the midpoint and Gauss-Legendre quadrature rules for a variety of different problems. For the midpoint rule, we found that the matrix $\mathbf{A}_{\text{S}}$ is of full rank for every problem as long as the restriction $N_{\text{eig}} \ge \text{max}(N_x,N_y) - 1$ is satisfied. However, for the Gauss-Legendre quadrature rule, we found that for some problems, the matrix $\mathbf{A}_{\text{S}}$ is rank deficient despite this \neweredit{inequality being satisfied}. We found that setting the value of $N_{\text{eig}}$ to be sufficiently large yields a matrix $\mathbf{A}_{\text{S}}$ of full rank. However, increasing the value of $N_{\text{eig}}$ is undesirable as the Gauss-Legendre rule \neweredit{can be} inaccurate for highly oscillatory functions. \finaledit{To summarise, \neweredit{on the basis of our findings we decided not to} use any form of Gaussian quadrature to approximate the integrals (\ref{eq:quadrature})--(\ref{eq:quadrature2}) as we cannot be sure that it will be accurate \neweredit{in our approach, bearing in mind that these approximations are key to the success of the semi-analytical method}. We therefore choose the midpoint rule to approximate the integrals (\ref{eq:quadrature})--(\ref{eq:quadrature2}). We note that more suitable quadrature methods may exist and/or could be developed, and if these methods can be implemented in the form of summations as presented in equations (\ref{eq:quadrature})--(\ref{eq:quadrature2}), then \finalrevision{they could be easily applied} to our semi-analytical method.}}
The midpoint rule uses evenly spaced abscissas, which leads to a restriction on the maximum number of terms, $N_{\text{eig}}$, taken in the summations appearing in the solution (\ref{eq:solution}). The functions {integrated in the definition of the coefficients} (\ref{eq:coeff_1_paper2})--(\ref{eq:coeff_4_paper2}) are similar to $\widetilde{f}(x) = f(x) \cos\left(k\pi x/L\right)$ with limits of integration of $x = 0$ and $x = L$. The period of $\widetilde{f}(x)$ is $2L/k$ and with $N_x$ evenly spaced abscissas used in the approximation of the integrals, the spacing between abscissas is $L/(N_x-1)$. If $2L/k = L/(N_x-1)$ the abscissas align with the peaks of $\widetilde{f}(x)$, which causes each evaluation of $\widetilde{f}(x)$ to occur at an abscissa $x = x^*$ at which $\widetilde{f}(x^*) = f(x^*)$.  This has the effect of returning the value of the integral {of $f(x)$ from $x = 0$ to $x = L$} instead of $\widetilde{f}(x)$. \finaledit{Hence, to ensure an accurate solution we set an upper limit on the maximum number of terms in the summations in equation (\ref{eq:solution}) to ensure $2L/k < L/(N_x-1)$, for all $k = 1,\hdots,N_{\text{eig}}$, leading to the restriction $N_{\text{eig}} \leq 2N_x-3$.  Combining both the lower and upper bounds on $N_{\text{eig}}$ yields the restriction:
\begin{align}
\label{eq:Neig_restriction}
\text{max}(N_x,N_y) - 1 \leq N_{\text{eig}} \leq 2\text{min}(N_x,N_y)-3.
\end{align}
We note that the lower bound in the restriction (\ref{eq:Neig_restriction}) applies to any quadrature rule, whereas the upper bound applies specifically to the midpoint rule.  \final{It also follows that $N_x \ge2$ must hold.} 
}



\section{Finite volume method}
\label{sec:fvm_paper2}
\absolutefinaledit{To verify the accuracy of our new semi-analytical method, we compare it {to} a standard vertex-centered finite volume method \citep{carr_2014}.} 
We first define a rectangular mesh over the \amend{unit cell $\mathcal{C}$} consisting of uniformly spaced nodes in the $x$ and $y$ directions. The nodes are located at $x = x_0 + kh_x$ for $k =0,\hdots,N_\text{F}^{(x)}$ and $y = y_0 + kh_y$ for $k =0,\hdots,N_\text{F}^{(y)}$, where $h_x = (x_n-x_0)/N_\text{F}^{(x)}$ and $h_y = (y_m-y_0)/N_\text{F}^{(y)}$. {Control volumes are formed around each node by connecting the centroid of each rectangular element to the midpoint of its edges. For interior nodes this gives a {rectangular control volume comprised of eight ``half-sized'' edges, each confined to a single element of homogeneous material.} The number of nodes {in the $x$ and $y$ directions} are always chosen to ensure that nodes coincide with interfaces between blocks so that each control volume edge is located within only one block.}

To {perform the finite volume discretization, the PDE (\ref{eq:diffusion}) is integrated over each control volume and the integrals along each control volume edge approximated using a midpoint rule. This procedure leads to} a {semi-discretized} finite volume equation of the form:
\begin{multline}
\label{eq:fvm}
\sum_{s \in e_x} D(x_{\text{mid}}^{(s)},y_{\text{mid}}^{(s)})\left(\frac{\partial \psi^{(\xi)}}{\partial x}(x_{\text{mid}}^{(s)},y_{\text{mid}}^{(s)}) + \delta_{x\xi} \right)\frac{h_x}{2} +\\ \sum_{s \in e_y} D(x_{\text{mid}}^{(s)},y_{\text{mid}}^{(s)})\left(\frac{\partial \psi^{(\xi)}}{\partial y}(x_{\text{mid}}^{(s)},y_{\text{mid}}^{(s)}) + \delta_{y\xi} \right)\frac{h_y}{2} = 0,
\end{multline}
where $e_x$ and $e_y$ are the sets of vertical and horizontal edges of the control volume, $(x_{\text{mid}}^{(s)},y_{\text{mid}}^{(s)})$ are the coordinates of the midpoint of edge $s$ of the control volume and $\delta_{x\xi}$ and $\delta_{y\xi}$ are Kronecker deltas, taking the value of one when the subscripts are equal and zero otherwise. {The equation} (\ref{eq:fvm}) is valid for all nodes, except those located along the boundaries $x = x_n$ and $y = y_m$. For these nodes, the finite volume equation (\ref{eq:fvm}) is modified to accommodate the periodic boundary conditions, as discussed in detail in \citep{carr_2013b}.

{Approximating the spatial derivatives appearing in equation (\ref{eq:fvm}) using second-order central differences and assembling the resulting finite volume equations for each node} yields a linear system of the form:
\begin{align}
\label{eq:linear_system_FVM}
\mathbf{A_{\text{F}}x_{\text{F}}} = \mathbf{b_{\text{F}}},
\end{align}
where $\mathbf{x_{\text{F}}}$ is a vector of length ${N_\text{F}^{(x)}N_\text{F}^{(y)}}$ containing the discrete unknown values of $\psi^{(\xi)}$ at each of the nodes. After solving this linear system, we use a bilinear interpolant of $\psi^{(\xi)}$ over each element to approximate the gradient within each element. This allows us to approximate the integrals {in equation (\ref{eq:deff_formula_1})--(\ref{eq:deff_formula_2})}, yielding the following formulae for the first and second columns of the {effective diffusivity}:
\begin{align}
\label{eq:deff_2}
[\mathbf{D}_{\mathrm{eff}}]_{(:,1)} = \frac{1}{A} \sum_{p= 1}^{N_\text{F}^{(x)}}\sum_{q = 1}^{N_\text{F}^{(y)}}D(x_{c}^{(p)},y_{c}^{(q)})(\nabla \psi^{(x)}(x_{c}^{(p)},y_{c}^{(q)})+\mathbf{e}_1)h_xh_y, \\ 
\label{eq:deff_3}
[\mathbf{D}_{\mathrm{eff}}]_{(:,2)} = \frac{1}{A} \sum_{p = 1}^{N_\text{F}^{(x)}}\sum_{q = 1}^{N_\text{F}^{(y)}}D(x_{c}^{(p)},y_{c}^{(q)})(\nabla \psi^{(y)}(x_{c}^{(p)},y_{c}^{(q)})+\mathbf{e}_2)h_xh_y,
\end{align}
where $(x_{c}^{(p)},y_{c}^{(q)})$ is the centroid of the $(p,q)$th element, $\mathbf{e}_1 = [1,0]^T$ and $\mathbf{e}_2 = [0,1]^T$.

\section{Results and \amend{discussion}}
\label{sec:results_paper2}
In this section we verify the accuracy and efficiency of our new semi-analytical method and demonstrate its applicability to computing effective diffusivities of complex {heterogeneous} geometries. In all comparisons between the semi-analytical and finite volume methods, we ensure that the spacing between abscissas and {the} spacing between nodes are identical {by setting $N_\text{F}^{(x)} = nN_x $ and $N_\text{F}^{(y)} = mN_y$ in the finite volume method, where $N_x$ and $N_y$ denote the number of abscissas per interface in the semi-analytical method}. Both the semi-analytical and finite volume methods were implemented in MATLAB and the linear systems for both methods{, equations (\ref{eq:linear_system_paper2}) and (\ref{eq:linear_system_FVM}), were stored in sparse format and solved using the backslash operator.\footnote{For both the semi-analytical and finite volume methods, the backslash operator implements the Unsymmetric MultiFrontal PACKage with automatic reordering.} All numerical experiments were carried out in MATLAB (version 2017A) running on an early 2015 MacBook Pro with a 2.7GHz dual-core Intel Core i5 processor and 8GB of memory.} 

The efficiency of the MATLAB code {implementing} the semi-analytical method can be {improved} by assuming that $m = n$ (i.e. the same number of blocks in both the horizontal and vertical directions) and that all $mn$ blocks are the same size. The main source of {the improvement} in efficiency is in being able to compute the summations appearing in the solution (\ref{eq:solution}) only once, instead of once for each abscissa used in the method. All of the {test cases} we consider and the majority of the problems we encountered in the literature have $m = n$ and identically sized blocks, so all calculations in this paper are computed using the {more efficient version of the} MATLAB code. The MATLAB codes implementing the semi-analytical method for both $m = n$ and $m \neq n$ are available on our GitHub repository: \href{https://github.com/NathanMarch/Homogenization}{https://github.com/NathanMarch/Homogenization}.

\subsection{Accuracy and efficiency of the semi-analytical method}
\label{sec:accuracy}
As a first verification of the accuracy of the semi-analytical solution we consider four geometries previously presented in \citep{szymkiewicz_2012}. These four test geometries are depicted in Table \ref{tab:effective_diffusivity} and take the form of a square cell of dimension $1$ by $1$ comprised of an $8$ by $8$ array of blocks each of dimension $0.125$ by $0.125$. The dark grey blocks have a diffusivity of $0.1$ and the light grey blocks have a diffusivity of $1$. Case 1 has a square shaped inclusion of dark grey blocks in the centre of the cell, case 2 consists of one larger dark grey inclusion in the centre of the cell and four smaller dark grey inclusions at the corners, case 3 consists of three layers, of which the middle layer is dark grey and case 4 contains an L-shaped dark grey inclusion. We compute the effective diffusivity using both the semi-analytical and finite volume methods for cases 1--4 and compare the results to those presented in \citep{szymkiewicz_2012}. We set $N_x = N_y = 16$ and {$N_\text{eig} = 2N_x - 3 = 29$ for the semi-analytical method and $N_{\text{F}}^{(x)} = N_{\text{F}}^{(y)} = 4N_{x} = 4N_{y} = 64$ for the finite volume method}. The results in Table \ref{tab:effective_diffusivity} demonstrate that both the semi-analytical and finite volume methods are in very good agreement with the results presented by \citep{szymkiewicz_2012}, agreeing to at least two significant figures for all entries in each {effective diffusivity}. 

\begin{table*}[hbt!]
\centering
\begin{tabular}{|l|c|c|c|}
\hline
 & Case 1 & Case 2 \\
\hline
& &\\[-0.3cm]
\multirow{-7}{*}{Geometry} & {\includegraphics[width=0.18\textwidth]{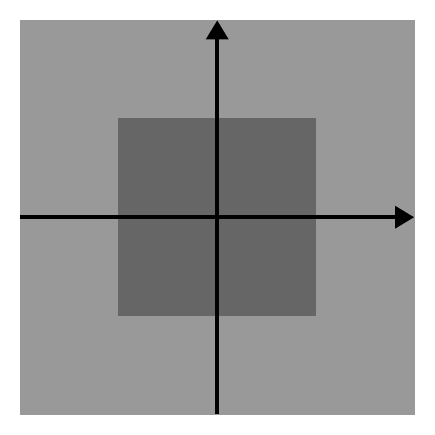}} & {\includegraphics[width=0.18\textwidth]{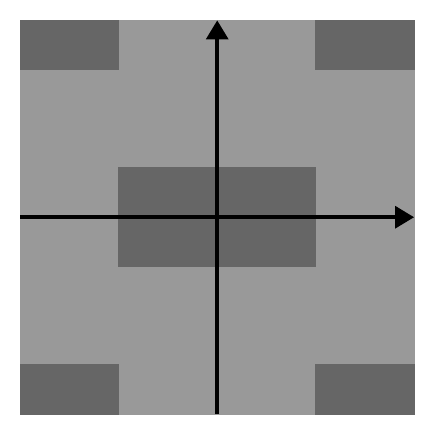}} \\
\hline
Szymkiewicz \citep{szymkiewicz_2012} & $ \begin{pmatrix} 0.649 &  0 \\ 0 & 0.649\end{pmatrix}$ & $ \begin{pmatrix} 0.693 &  0 \\ 0 & 0.566\end{pmatrix}$\\ 
Finite Volume & $ \begin{pmatrix} 0.648 &  0 \\ 0 & 0.648\end{pmatrix}$ & $ \begin{pmatrix} 0.694 &  0 \\ 0 & 0.566\end{pmatrix}$\\
Semi-Analytical & $ \begin{pmatrix} 0.648 &  0 \\ 0 & 0.648\end{pmatrix}$ & $ \begin{pmatrix} 0.693 &  0 \\ 0 & 0.566\end{pmatrix}$\\
\hline
 & Case 3 & Case 4\\
 \hline
 & &\\[-0.3cm]
\multirow{-7}{*}{Geometry}  & {\includegraphics[width=0.18\textwidth]{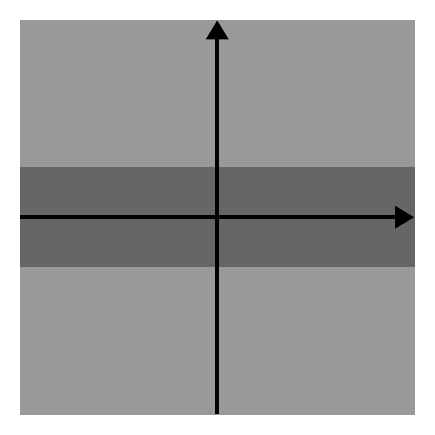}} & {\includegraphics[width=0.18\textwidth]{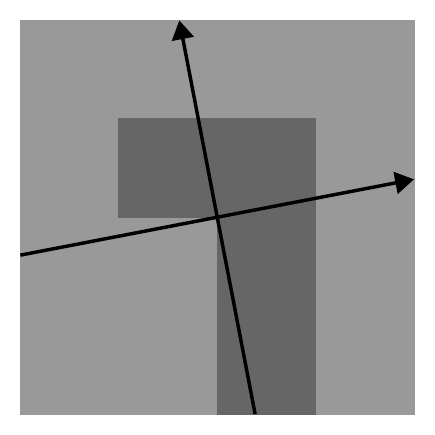}}\\
\hline
Szymkiewicz \citep{szymkiewicz_2012}  & $ \begin{pmatrix} 0.775 &  0 \\ 0 & 0.309\end{pmatrix}$ & $ \begin{pmatrix} 0.533 &  -0.0286 \\ -0.0286 & 0.675 \end{pmatrix}$ \\ 
Finite Volume & $ \begin{pmatrix} 0.775 &  0 \\ 0 & 0.308\end{pmatrix}$ & $ \begin{pmatrix} 0.533 &  -0.0286 \\ -0.0286 & 0.676\end{pmatrix}$\\
Semi-Analytical  &$ \begin{pmatrix} 0.775 &  0 \\ 0 & 0.308\end{pmatrix}$ & $ \begin{pmatrix} 0.533 &  -0.0286 \\ -0.0286 & \revision{0.676}\end{pmatrix}$\\
\hline
\end{tabular}
\caption{Effective diffusivity for cases 1--4 as calculated using the semi-analytical and finite volume methods and compared to the results presented in \citep{szymkiewicz_2012}. The principal directions of diffusion are overlaid upon each geometry. All results are {reported} to three significant figures, {consistent with} \citep{szymkiewicz_2012}.}
\label{tab:effective_diffusivity}
\end{table*}

\begin{table*}[!htb]
\centering
\def\arraystretch{1.05}
\begin{tabular}{|c|c|c|c|c|}
\hline
Case A & \multicolumn{2}{c}{Semi-Analytical} \vline & \multicolumn{2}{c} {Finite Volume} \vline \\
\hline
$N_x$ & Relative Error ($\mathbf{E}$) & Time (s) & Relative Error ($\mathbf{E}$) &Time (s) \\
\hline
4 & $ \begin{pmatrix}  4.09 \text{e}{-03} &  0 \\ 0 & 4.09 \text{e}{-03} \end{pmatrix}$ & 0.00805 & $ \begin{pmatrix} 7.85 \text{e}{-03} &  0 \\ 0 & 7.85 \text{e}{-03} \end{pmatrix}$ & 0.00938\\
8 & $ \begin{pmatrix} 1.80 \text{e}{-03} &  0 \\ 0 & 1.80 \text{e}{-03}\end{pmatrix}$ & 0.0112 &$ \begin{pmatrix} 2.89 \text{e}{-03} &  0 \\ 0 & 2.89 \text{e}{-03}\end{pmatrix}$ & 0.0280\\ 
16 & $ \begin{pmatrix} 8.34 \text{e}{-04} &  0 \\ 0 & 8.34 \text{e}{-04}\end{pmatrix}$ & 0.0338 & $ \begin{pmatrix} 1.05 \text{e}{-03} &  0 \\ 0 & 1.05 \text{e}{-03}\end{pmatrix}$ & 0.113\\ 
32 & $ \begin{pmatrix} 4.04 \text{e}{-04} &  0 \\ 0 & 4.04 \text{e}{-04}\end{pmatrix}$ & 0.0630 &$ \begin{pmatrix} 3.69 \text{e}{-04} &  0 \\ 0 & 3.69 \text{e}{-04}\end{pmatrix}$ & 0.531\\
64 & $ \begin{pmatrix} 2.04  \text{e}{-04} &  0 \\ 0 & 2.04 \text{e}{-04}\end{pmatrix}$ & 0.269 & $ \begin{pmatrix} 1.22 \text{e}{-04} &  0 \\ 0 & 1.22 \text{e}{-04}\end{pmatrix}$ & 2.93\\ 
\hline
Case B & \multicolumn{2}{c}{Semi-Analytical}  \vline & \multicolumn{2}{c} {Finite Volume} \vline \\
\hline
$N_x$ & Relative Error ($\mathbf{E}$) & Time (s) & Relative Error ($\mathbf{E}$) &Time (s) \\
\hline
4 & $ \begin{pmatrix} 6.84 \text{e}{-03} &  5.04 \text{e}{-03} \\ 5.04 \text{e}{-03} & 4.47 \text{e}{-03}\end{pmatrix}$ & 0.00747 & $ \begin{pmatrix} 1.30 \text{e}{-02} &  2.44 \text{e}{-03} \\ 2.44 \text{e}{-03} & 8.47 \text{e}{-03}\end{pmatrix}$ & 0.00923\\ 
8 & $ \begin{pmatrix} 3.01 \text{e}{-03} &  2.21 \text{e}{-03} \\ 2.21 \text{e}{-03} & 1.98 \text{e}{-03}\end{pmatrix}$ & 0.0109 & $\begin{pmatrix} 4.82 \text{e}{-03} &  1.88 \text{e}{-03} \\ 1.88 \text{e}{-03} & 3.14 \text{e}{-03}\end{pmatrix}$ & 0.0277\\ 
16 &$ \begin{pmatrix} 1.40 \text{e}{-03} &  1.02 \text{e}{-03} \\ 1.02 \text{e}{-03} & 9.23 \text{e}{-04}\end{pmatrix}$ & 0.0331 & $ \begin{pmatrix} 1.75 \text{e}{-03} &  9.12 \text{e}{-04} \\ 9.12 \text{e}{-04} & 1.14 \text{e}{-03}\end{pmatrix}$ & 0.115\\ 
32 & $ \begin{pmatrix} 6.77 \text{e}{-04} &  4.94 \text{e}{-04} \\ 4.94 \text{e}{-04} & 4.48 \text{e}{-04}\end{pmatrix}$ & 0.0629 & $ \begin{pmatrix} 6.17 \text{e}{-04} &  3.76 \text{e}{-04} \\ 3.76 \text{e}{-04} & 4.02 \text{e}{-04}\end{pmatrix}$ & 0.530\\
64 & $ \begin{pmatrix} 3.42 \text{e}{-04} & 2.50 \text{e}{-04} \\ 2.50 \text{e}{-04} & 2.27 \text{e}{-04}\end{pmatrix}$ & 0.270 & $ \begin{pmatrix} 2.05 \text{e}{-04} &  1.36 \text{e}{-04} \\ 1.36 \text{e}{-04} & 1.33 \text{e}{-04}\end{pmatrix}$ & 2.92\\ 
\hline
Case C & \multicolumn{2}{c}{Semi-Analytical}  \vline & \multicolumn{2}{c} {Finite Volume} \vline \\
\hline
$N_x$ & Relative Error ($\mathbf{E}$) & Time (s) & Relative Error ($\mathbf{E}$) &Time (s) \\
\hline
4 & $ \begin{pmatrix} 5.60 \text{e}{-04} &  0 \\ 0 & 5.60 \text{e}{-04} \end{pmatrix}$ & 0.0702& $ \begin{pmatrix} 1.05 \text{e}{-03} &  0 \\ 0 & 1.05 \text{e}{-03} \end{pmatrix}$ & 0.124\\ 
8 & $ \begin{pmatrix} 2.52 \text{e}{-04} &  0 \\ 0 & 2.52 \text{e}{-04}\end{pmatrix}$ & 0.151 &$ \begin{pmatrix} 3.69 \text{e}{-04} &  0 \\ 0 & 3.69 \text{e}{-04}\end{pmatrix}$ & 0.549\\ 
16 & $ \begin{pmatrix} 1.25 \text{e}{-04} &  0 \\ 0 & 1.25 \text{e}{-04}\end{pmatrix}$ & 0.619 & $ \begin{pmatrix} 1.22 \text{e}{-04} &  0 \\ 0 & 1.22 \text{e}{-04}\end{pmatrix}$ & 2.98\\ 
\hline
Case D & \multicolumn{2}{c}{Semi-Analytical} \vline & \multicolumn{2}{c} {Finite Volume} \vline \\
\hline
$N_x$ & Relative Error ($\mathbf{E}$) & Time (s) & Relative Error ($\mathbf{E}$) &Time (s) \\
\hline
4 & $ \begin{pmatrix} 9.38 \text{e}{-04} &  7.42 \text{e}{-04} \\ 7.42 \text{e}{-04} & 6.10 \text{e}{-04}\end{pmatrix}$ & 0.0696 & $ \begin{pmatrix} 1.75 \text{e}{-03} &  9.12 \text{e}{-04} \\ 9.12 \text{e}{-04} & 1.14 \text{e}{-03}\end{pmatrix}$ & 0.0650\\ 
8 & $ \begin{pmatrix} 4.22 \text{e}{-04} &  3.34 \text{e}{-04} \\ 3.34 \text{e}{-04} & 2.75 \text{e}{-04}\end{pmatrix}$ & 0.151 & $\begin{pmatrix} 6.17 \text{e}{-04} &  3.76 \text{e}{-04} \\ 3.76 \text{e}{-04} & 4.02 \text{e}{-04}\end{pmatrix}$ & 0.388\\ 
16 &$ \begin{pmatrix} 2.10 \text{e}{-04} &  1.65 \text{e}{-04} \\ 1.65 \text{e}{-04} & 1.36 \text{e}{-04}\end{pmatrix}$ & 0.615 & $ \begin{pmatrix} 2.05 \text{e}{-04} &  1.36 \text{e}{-04} \\ 1.36 \text{e}{-04} & 1.33 \text{e}{-04}\end{pmatrix}$ & 2.43\\ 
\hline
\end{tabular}
\caption{Error {and runtimes} of the {effective diffusivities} calculated using the semi-analytical and finite volume methods.} 
\label{tab:case1_4}
\end{table*}

The {effective diffusivities} for cases 1--3 are diagonal matrices as each geometry is invariant under a rotation of $180^{\circ}$ \citep{sviercoski_2010}. Additionally, the diagonal entries of the {effective diffusivity} for case 1 are identical as the geometry is {also} invariant under a rotation of $90^{\circ}$. As cases 1--3 have {effective diffusivities} that are diagonal, the principal directions of diffusion are in the horizontal ($x)$ and vertical ($y$) directions and the diagonal elements represent the diffusivity in these directions {(see Table \ref{tab:effective_diffusivity})}. For cases 2 and 3 the diffusivity in the horizontal direction is larger than the diffusivity in the vertical direction. This is because, unlike for flow in the horizontal direction, flow in the vertical direction must pass through or around the low diffusivity (dark grey) blocks. Case 3 is a layered medium, so the diffusivity in the horizontal and vertical directions {are the area-weighted arithmetic and harmonic averages of the layer diffusivities, respectively}. The eigenvalues of the {effective diffusivity} for case 4 are $0.527$ and $0.682$ and the corresponding normalised eigenvectors are $[0.982,0.189]^T$ and $[0.189,-0.982]^T$, {which define principal directions of diffusion for case 4 that are rotated} $10.9^{\circ}$ anti-clockwise from the standard Cartesian axes {(see Table \ref{tab:effective_diffusivity})}.

We now compare the accuracy of the semi-analytical {method} to the finite volume {method} by first computing a benchmark {effective diffusivity} using the finite volume method with a very fine node spacing. We consider the geometries of cases 1 and 4 as presented in \citep{szymkiewicz_2012}, represented as both $4$ by $4$ and $16$ by $16$ grids of {equal-sized} blocks. \final{We label theses cases as cases A--D such that:
\finalrevision{\begin{itemize}
\item Case A is the $4$ by $4$ configuration of Case 1;
\item Case B is the $4$ by $4$ configuration of Case 4;
\item Case C is the $16$ by $16$ configuration of Case 1;
\item Case D is the $16$ by $16$ configuration of Case 4.
\end{itemize}}

 {For cases A and B each block is of size $0.25$ by $0.25$ while for the $16$ by $16$ configuration each block is of size $0.0625$ by $0.0625$. For cases A and B we compute the effective diffusivity using $N_x = N_y = 4, 8, 16, 32$ and $64$ for the semi-analytical method and $N_{\text{F}}^{(x)} = N_{\text{F}}^{(y)} = 4N_x = 4N_y = 16, 32, 64, 128$ and $256$ for the finite volume method. For cases C and D we compute the effective diffusivity using $N_x = N_y = 4, 8$ and $16$ for the semi-analytical method and $N_{\text{F}}^{(x)} = N_{\text{F}}^{(y)} = 16N_x = 16N_{y} = 64, 128$ and $256$ for the finite volume method. For all cases}, the benchmark effective diffusivity is computed using the finite volume method with $N_{\text{F}}^{(x)} = N_{\text{F}}^{(y)} = 1024$}. For the semi-analytical method, {$N_{\text{eig}}$} is set {as} the minimum of $2N_x-3$ and $100$. This is because, as mentioned in section \ref{sec:implementation}, {the value} $2N_x-3$ is {the upper limit on the maximum number of terms in the summations in equation (\ref{eq:solution})} and in preliminary testing we found that using more than $100$ terms in the summation had negligible effect on the solution. We compute our relative error $\mathbf{E} =  |(\mathbf{D}_{\mathrm{eff}}  -\widehat{\mathbf{D}}_{\text{eff}}) \varoslash  \widehat{\mathbf{D}}_{\text{eff}}|$, where $\varoslash$ represents Hadamard/element-wise division,  $\mathbf{D}_{\mathrm{eff}}$ represents the effective diffusivity as calculated from the semi-analytical/finite volume solution under comparison and $\widehat{\mathbf{D}}_{\text{eff}}$ is the benchmark effective diffusivity. We also record the runtime for both the semi-analytical and finite volume methods, which includes the time taken to solve the boundary value problem {(\ref{eq:diffusion_3})--(\ref{eq:zero_mean_2})} and compute the effective diffusivity via equations (\ref{eq:eff_diff_1})--(\ref{eq:eff_diff_4}) and (\ref{eq:deff_2})--(\ref{eq:deff_3}). All runtimes reported in this section were calculated by performing each calculation ten times and taking the median runtime. {Results in Table \ref{tab:case1_4} demonstrate that} the semi-analytical method is either faster or more accurate (or both) than the finite volume method, indicating that for a specified level of accuracy the semi-analytical method is more efficient. 

\begin{figure}[p]
\centering
\includegraphics[width=0.43\textwidth]{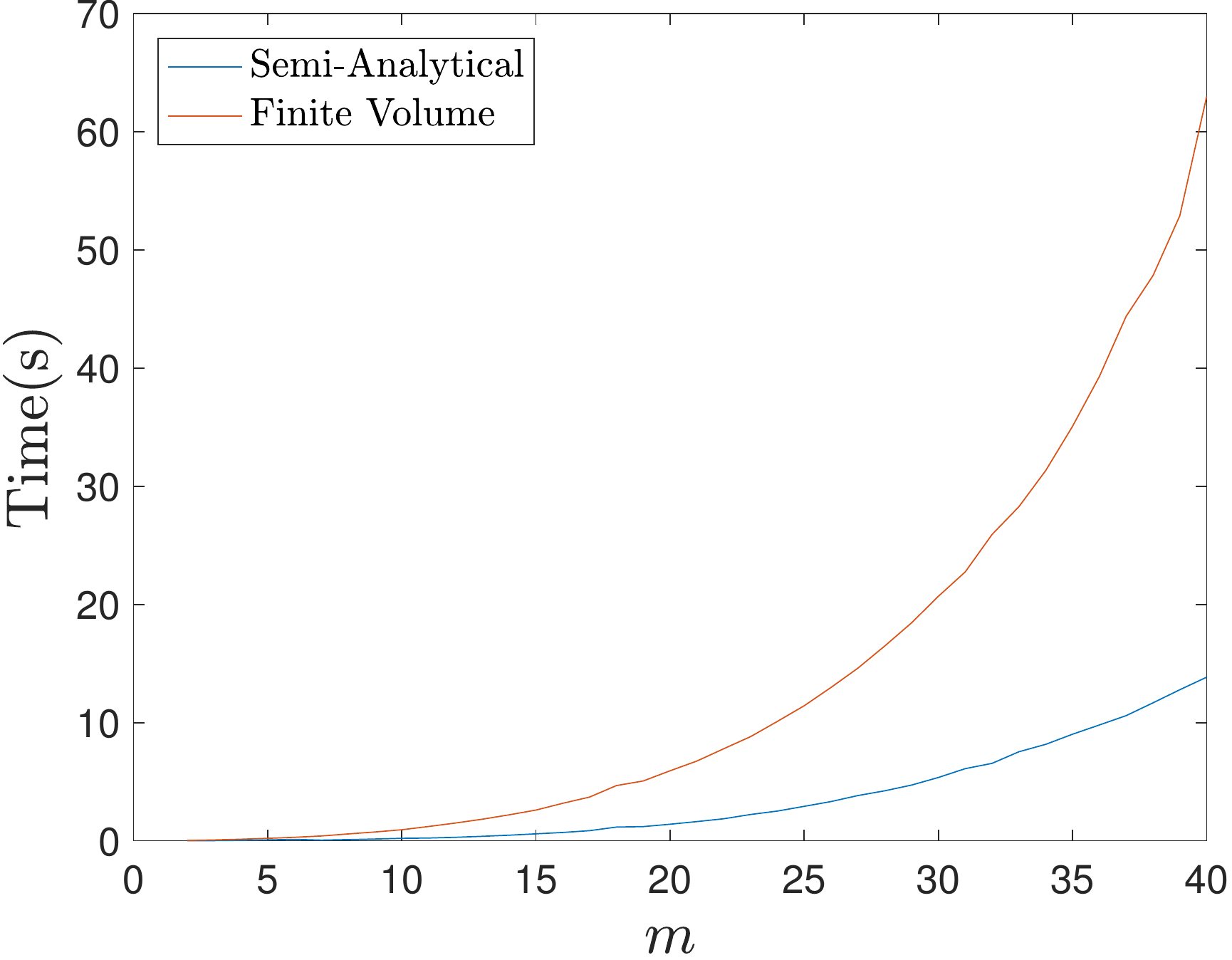}
\includegraphics[width=0.43\textwidth]{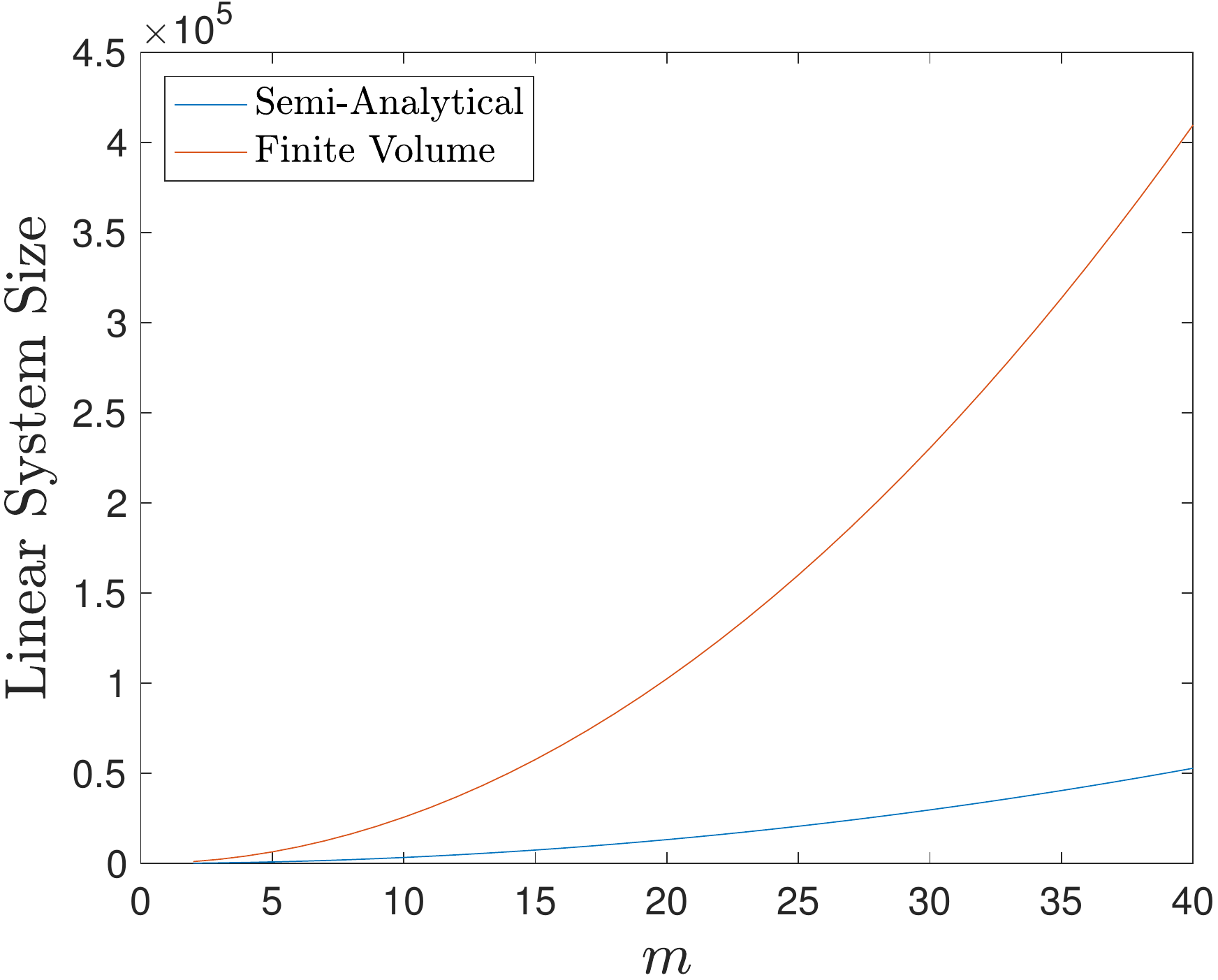}
\caption{Comparison of time taken to solve the linear system and size of the linear system for both the semi-analytical and finite volume methods for geometries consisting of an $m$ by $m$ grid, for $m = 2,\hdots,40$. The time taken represents the time to formulate and solve the linear systems, as well as to calculate the effective diffusivity.}
\label{fig:size_times}
\end{figure}

\begin{table}[p] 
\centering
\begin{tabular}{|l|c|c|c|}
\hline
 & $10 \times 10$ & $20 \times 20$  \\
\hline
& &\\[-0.3cm]
\multirow{-10}{*}{\edit{Geometry}} & {\includegraphics[width=0.26\textwidth,height = 0.26\textwidth]{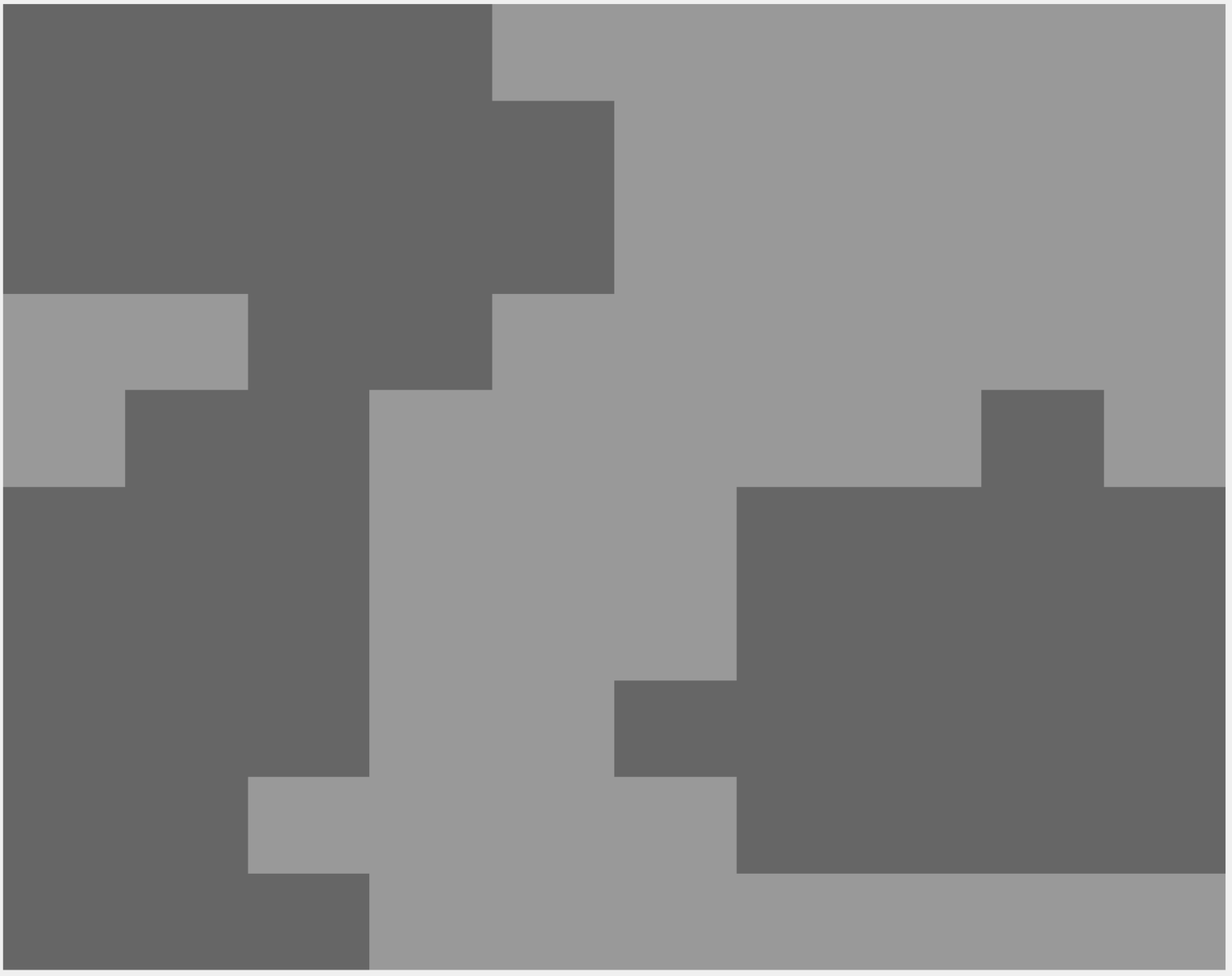}} & {\includegraphics[width=0.26\textwidth,height = 0.26\textwidth]{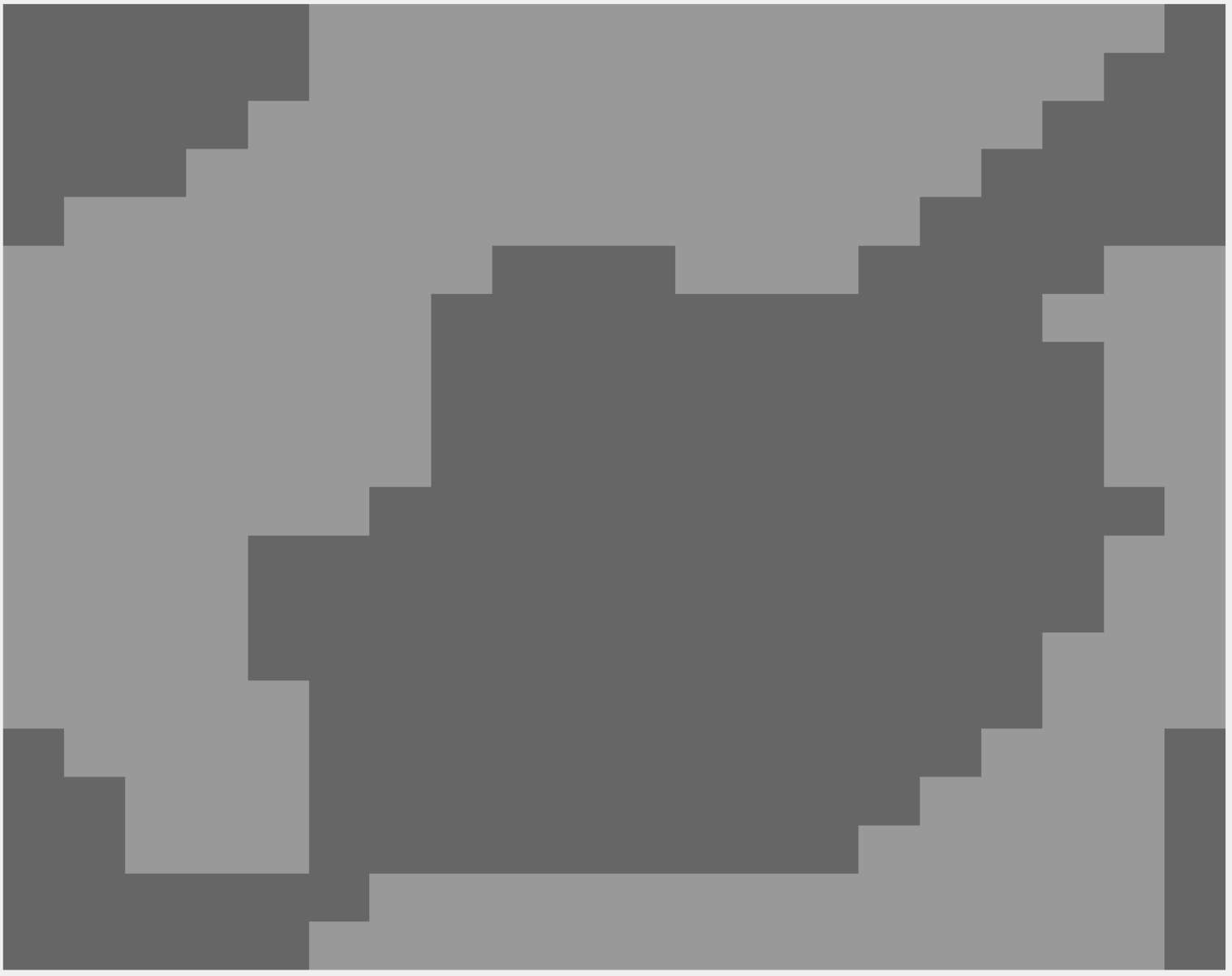}}\\
\hline
Effective diffusivity& $ \begin{pmatrix} 0.272 &  0.00754 \\ 0.00754 & 0.374\end{pmatrix}$ & $ \begin{pmatrix} 0.341 &  -0.0920 \\ -0.0920 & 0.346 \end{pmatrix}$\\ 
\hline
\hline
 & $50 \times 50$  & $100 \times 100$  \\
\hline
& &\\[-0.3cm]
\multirow{-10}{*}{\edit{Geometry}} & {\includegraphics[width=0.26\textwidth,height = 0.26\textwidth]{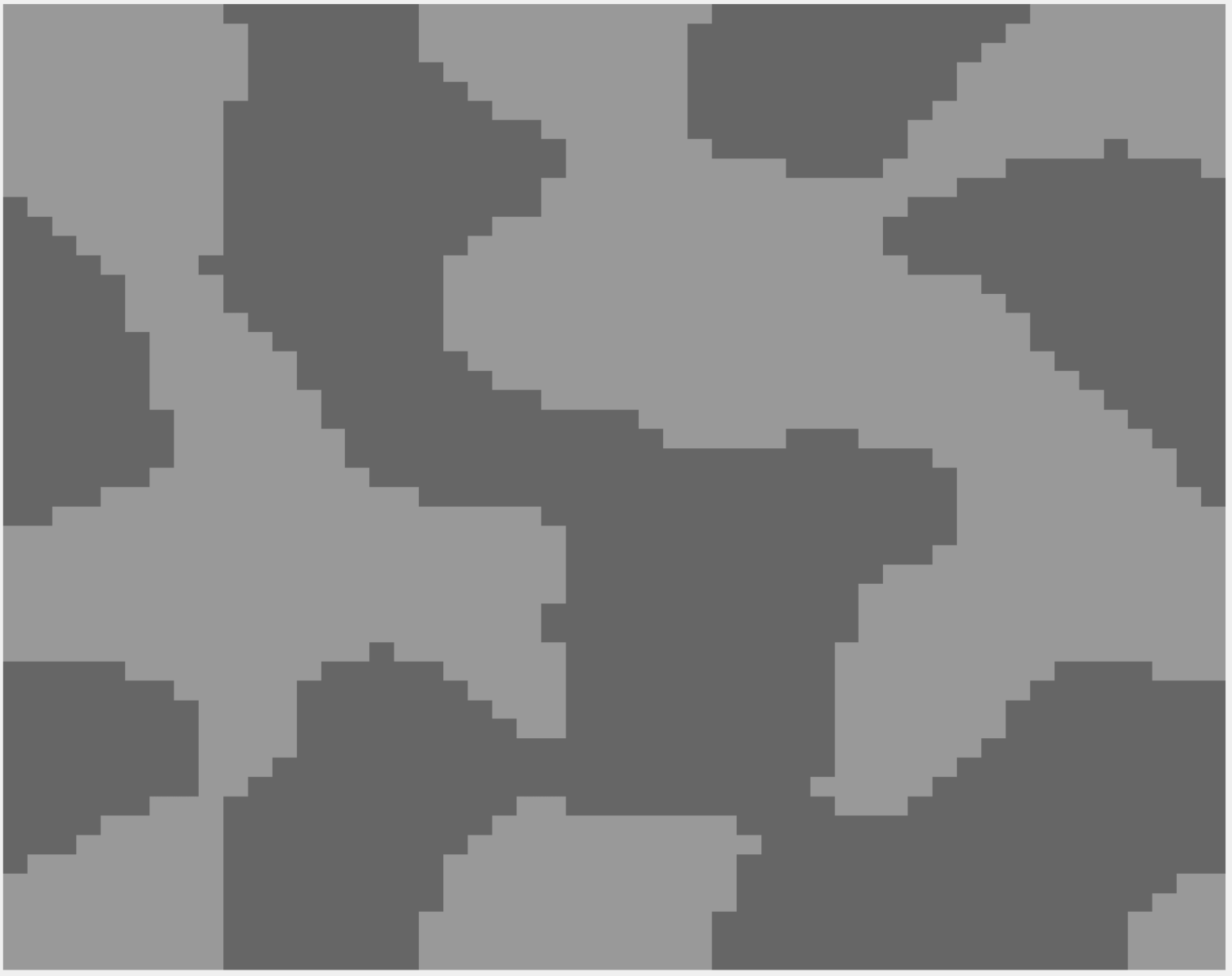}} & {\includegraphics[width=0.26\textwidth,height = 0.26\textwidth]{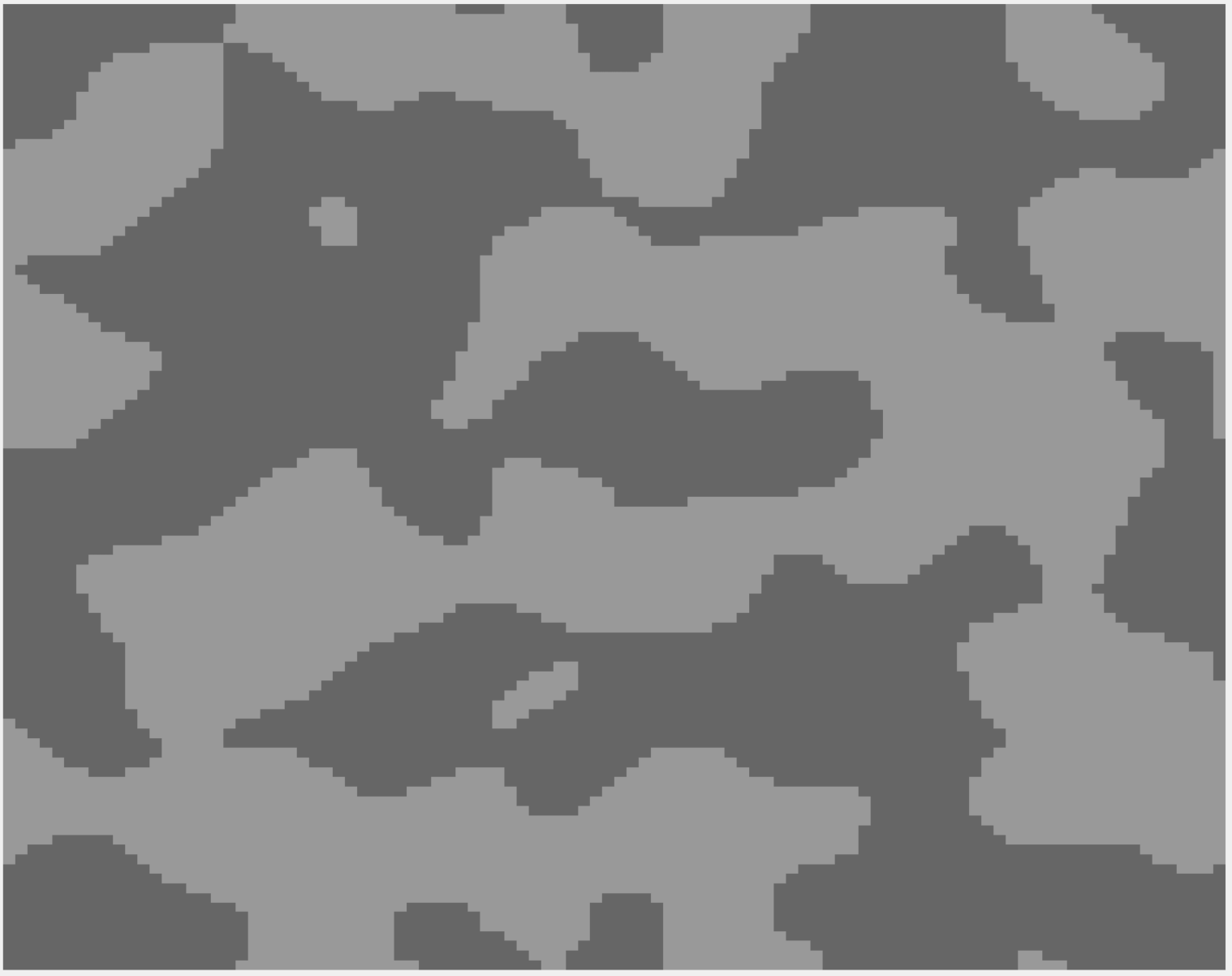}} \\
\hline
Effective diffusivity & $ \begin{pmatrix} 0.310&  0.0177 \\ 0.0177 & 0.342\end{pmatrix}$  & $ \begin{pmatrix} 0.340 &  0.000954  \\ 0.000954 & 0.304 \end{pmatrix}$\\ 
\hline
\end{tabular}
\caption{Effective diffusivity for {four} randomly generated {aggregated media consisting of a $10$ by $10$, $20$ by $20$, $50$ by $50$ and $100$ by $100$ array of blocks, respectively,} as calculated using the semi-analytical method. In all four cases {$N_x = N_y = 8$ and} the number of blocks with diffusivity 0.1 (dark grey blocks) is equal to the number of blocks with diffusivity 1 (light grey blocks). The aggregated geometries were inspired by those presented in \citep{ray_2017}.}
\label{tab:random_geometries}
\end{table}

We now further explore the efficiency of our new semi-analytical method by comparing the size of the linear systems, equations (\ref{eq:linear_system_paper2}) and (\ref{eq:linear_system_FVM}), and the time taken to compute the effective diffusivity for both the semi-analytical and finite volume methods. To do this we consider a square \revision{unit cell}  with $x_0 = y_0 = 0$ and $x_n = y_m = 1$, consisting of an $m$ by $m$ grid of {equal-sized} blocks for $m = 2,\hdots,40$, i.e.{,} we consider a $2$ by $2$ grid, a $3$ by $3$ grid etc. The \revision{unit cell} is a checkerboard structure consisting of diffusivities of $0.1$ and $1$ with a diffusivity of $0.1$ in the top left corner. For example, the diffusivities for $m = 3$ are:
\begin{gather}
\label{eq:checkerboard}
\begin{bmatrix}
D_{1,1} & D_{1,2} & D_{1,3}\\
D_{2,1} & D_{2,2} & D_{2,3}\\
D_{3,1} & D_{3,2} & D_{3,3}\\
\end{bmatrix} = \begin{bmatrix}
0.1 & 1 & 0.1\\
1 & 0.1 & 1\\
0.1 & 1 & 0.1\\
\end{bmatrix}.
\end{gather}
As usual, we specify the spacing between abscissas for the semi-analytical method to be the same as the spacing between nodes for the finite volume method. {We compute the effective diffusivity using $N_x = N_y = 16$ and $N_\text{eig} = 2N_x - 3 = 29$ for the semi-analytical method and $N_{\text{F}}^{(x)} = N_{\text{F}}^{(y)} = 16m$ for the finite volume method. \finaledit{This leads to the linear system (\ref{eq:linear_system_paper2}) in the semi-analytical method having dimension $m^{2}(2N_{x}+1) = 33m^{2}$ and the linear system (\ref{eq:linear_system_FVM}) in the finite volume method having dimension $m^{2}N_{x}^{2} = 256m^{2}$}. Results in Figure \ref{fig:size_times} demonstrate that the efficiency advantage of the semi-analytical method over the finite volume method becomes more pronounced as $m$ increases. This is because as the number of blocks increases, the size of the linear system increases (see Figure \ref{fig:size_times}) and the time taken to solve these systems becomes the dominant computational cost {of calculating the effective diffusivity}.} This {improvement} in efficiency that the semi-analytical method has over the finite volume method demonstrates the potential of the semi-analytical method to significantly speed up simulations of transport processes across block, heterogeneous media.

\finaledit{\subsection{Accuracy of the solution of the linear system (\ref{eq:linear_system_paper2})}
In this section we demonstrate the accuracy of the solution of the linear system (\ref{eq:linear_system_paper2}) corresponding to cases 1--4 with $N_x = N_y = 16$ and all values of $N_{\text{eig}}$ that satisfy the restriction (\ref{eq:Neig_restriction}), which are $N_{\text{eig}} = 15,\hdots,29$. We stress that it is important that the lower bound of the restriction (\ref{eq:Neig_restriction}) is satisfied, as otherwise the coefficient matrix $\mathbf{A}_s$ will be singular. For each case, the condition number (as calculated using the $2$-norm) of the coefficient matrix $\mathbf{A}_s$ is of the order $10^3$ and the $2$-norm of the residual of the linear system (\ref{eq:linear_system_paper2}) is of the order $10^{-14}$, which means the continuity conditions (\ref{eq:ICX1})--(\ref{eq:ICX4}) are satisfied with a maximum absolute error of $10^{-14}$. An upper bound on the relative error in the solution of the linear system (\ref{eq:linear_system_paper2}) is $\text{cond}(\mathbf{A}_s)\norm{\mathbf{r}_s}_2/\norm{\mathbf{b}_s}_2$, where $\mathbf{r}_s$ is the residual of the linear system (\ref{eq:linear_system_paper2}), which yields a maximum relative error of \finalfinaledit{$10^{-11}/\norm{\mathbf{b}_s}_2$.}}

\subsection{Results for \change{complex structured heterogeneous} geometries}
\label{sec:complex_geometry}
\finaledit{We now demonstrate how the semi-analytical method can be used to homogenize complex heterogeneous geometries.} For our demonstration, we consider a square cell of dimension $1$ by $1$ consisting of an $m$ by $m$ array of square blocks and generate aggregated, random geometries of two different media. We first {generate $m^{2}$ uniformly distributed, random numbers between $0$ and $1$ denoted by $D^{(0)}_{i,j}$ where $i = 1,\hdots, m$ and $j = 1,\hdots,m$\finaledit{.} We then perform the following iteration:
\begin{multline}
\label{eq:iteration}
D^{(k+1)}_{i,j} = \frac{4}{9} D^{(k)}_{i,j} + \frac{1}{9} \left(D^{(k)}_{i-1,j} +   D^{(k)}_{i+1,j} +   D^{(k)}_{i,j-1} +   D^{(k)}_{i,j+1}\right)  \\+ \frac{1}{36}\left(D^{(k)}_{i-1,j-1} + D^{(k)}_{i-1,j+1} + D^{(k)}_{i+1,j-1} + D^{(k)}_{i+1,j+1}\right),
\end{multline}
for all $i = 1,\hdots, m$ and $j = 1,\hdots,m$, where $D_{i,j}^{(k)}$} represents the diffusivity in the $(i,j)$th block of the \revision{unit cell $\mathcal{C}$} (see Figure \ref{fig:figure1}(b)) at the $k$th iteration and {the weightings are inspired by those used in lattice Boltzmann methods \citep{perumal_2015}.} For the purposes of the algorithm in equation (\ref{eq:iteration}), {we assume periodicity to process} the blocks along the boundaries{, for example,} when $i = 1$, we set $i-1$ to be equal to $m$. After performing a fixed number of iterations, we are left with an array of diffusivities containing entries between $0$ and $1$, where larger values are aggregated together and smaller values are aggregated together. \finaledit{The reason for this aggregation is that during each iteration, the diffusivity of a given block is a weighted average of itself and its nearest neighbours, so in each iteration, the diffusivities are converging to those of its nearest neighbours.} {To ensure the geometry is comprised of equal parts light grey (diffusivity of 1) and dark grey (diffusivity of 0.1), we then prescribe that the largest $0.5 m^{2}$ entries of the array are 1 and the remaining $0.5m^{2}$ entries are 0.1. 

Table \ref{tab:random_geometries} presents some example geometries generated using the above procedure for $m = 10, 20, 50$ and $100$, which resemble realistic geometries presented in \citep{amaziane_2001} and \citep{sviercoski_2010}. These figures demonstrate the highly complex heterogeneous geometries that can be captured using blocks and how the semi-analytical method can be used to calculate the effective diffusivity for such geometries.} \finaledit{We note that each iteration in the aggregation procedure could be considered as a change in the underlying geometry of a domain \finalfinaledit{arising} throughout the simulation of a diffusive transport problem. In this scenario, the effective diffusivity can be calculated easily and efficiently using our semi-analytical method, as it requires merely changing the diffusivities in some of the blocks. We also note that our semi-analytical method can be applied to any geometry where the diffusivities are represented as a grid of pixels (e.g. CT scans).} 

\subsection{\change{Results for complex unstructured heterogeneous geometries}}
\label{sec:pixelation}
We now demonstrate the applicability of the semi-analytical method to geometries involving curved \change{edges} by first pixelating them, which allows for the geometry to be represented as a grid of blocks. \revision{We consider a geometry presented in \citep{carr_2014} (see Table \amend{\ref{tab:pixelation_refinement}(a)}) and pixelate it several times to generate geometries of different sizes (see Figure \ref{fig:figure1} (b)). This is a $512$ by $512$ geometry in which the dark grey blocks have a diffusivity of $0.1$ and the light grey blocks have a diffusivity of $1$.} \finaledit{We define a pixelation parameter $r$, such that the image is now represented as an $r$ by $r$ grid of square blocks, so decreasing values of $r$ \revision{produce greater pixelation}. Each block in the pixelated geometry is made up of a $k$ by $k$ grid of blocks, such that $rk = m$ \revision{(for the geometry in \amend{Figure \ref{fig:figure1}(a)}, $m = 512$)}. On each $k$ by $k$ grid of blocks, we compute the mean of the diffusivities \revision{of} the $k^2$ blocks\revision{. If} the mean is greater than or \amend{equal to $0.55$}, the diffusivity of the \revision{entire $k$ by $k$ block} is set to be $1$, otherwise it is set to be $0.1$. The threshold value of $0.55$ is chosen as it is halfway between the low diffusivity of $0.1$ and the high diffusivity of $1$ and as the pixellation is performed via a simple averaging process, the computational cost is negligible.} On each of the pixelated geometries, we compute the effective diffusivity  using the semi-analytical method with $N_x = N_y = 4$ \revision{abscissas per interface} and \finaledit{compare the result to a benchmark effective diffusivity obtained from the finite volume method of section  \ref{sec:fvm_paper2} applied to the original 512 by 512 geometry with a fine mesh \finaledit{($N^{(x)}_{\text{F}}$ = $N^{(y)}_{\text{F}}$ = \final{$1024$})}.} We consider levels of pixelation corresponding to \amend{$r = 16,32,64,128$} and for each level of pixelation compute the relative error $\mathbf{E}$ \revision{defined in section \ref{sec:accuracy}}.

{{Results in Table \ref{tab:pixelation_refinement} demonstrate} that as the pixelation increases (as $r$ decreases), the error in the effective {diffusivity} \amend{generally increases, that is, higher values of $r$ generally yield more accurate effective diffusivities.}  These results demonstrate the applicability of pixelating complicated geometries to the fast computation of accurate effective {diffusivities}. {For example}, pixelating the \amend{image} to a $32$ by $32$ \amend{geometry} {allows} the \amend{diagonal entries of the} effective {diffusivity} to be calculated using our semi-analytical method to an accuracy of approximately $1\%$ in under $1$ second. {Although a standard finite volume method could {also} be applied to the pixelated $32$ by $32$ geometry, our {semi-analytical method} is faster and {of equivalent accuracy}, as demonstrated in Table \ref{tab:case1_4}.}} 

\subsection{Extension of semi-analytical method to geometries without aligned interfaces}
\label{sec:esa}
\begin{figure}[htbp!]
\centering 
\subfloat[\revision{Semi-analytical method abscissas}]{\includegraphics[width=0.43\textwidth,height=0.43\textwidth]{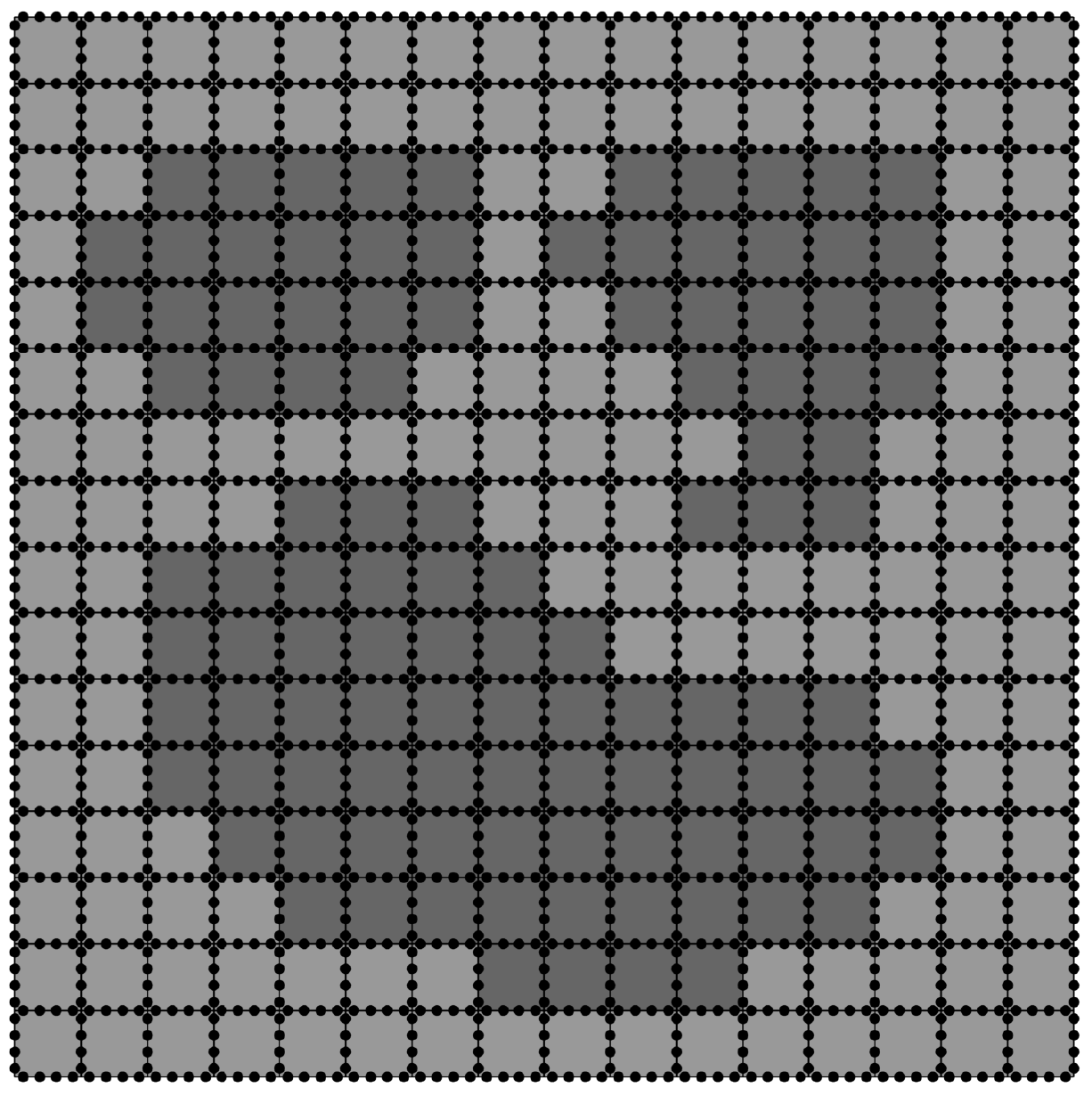}}\hspace{0.2cm}
\subfloat[Extended semi-analytical method abscissas]{\includegraphics[width=0.43\textwidth,height=0.43\textwidth]{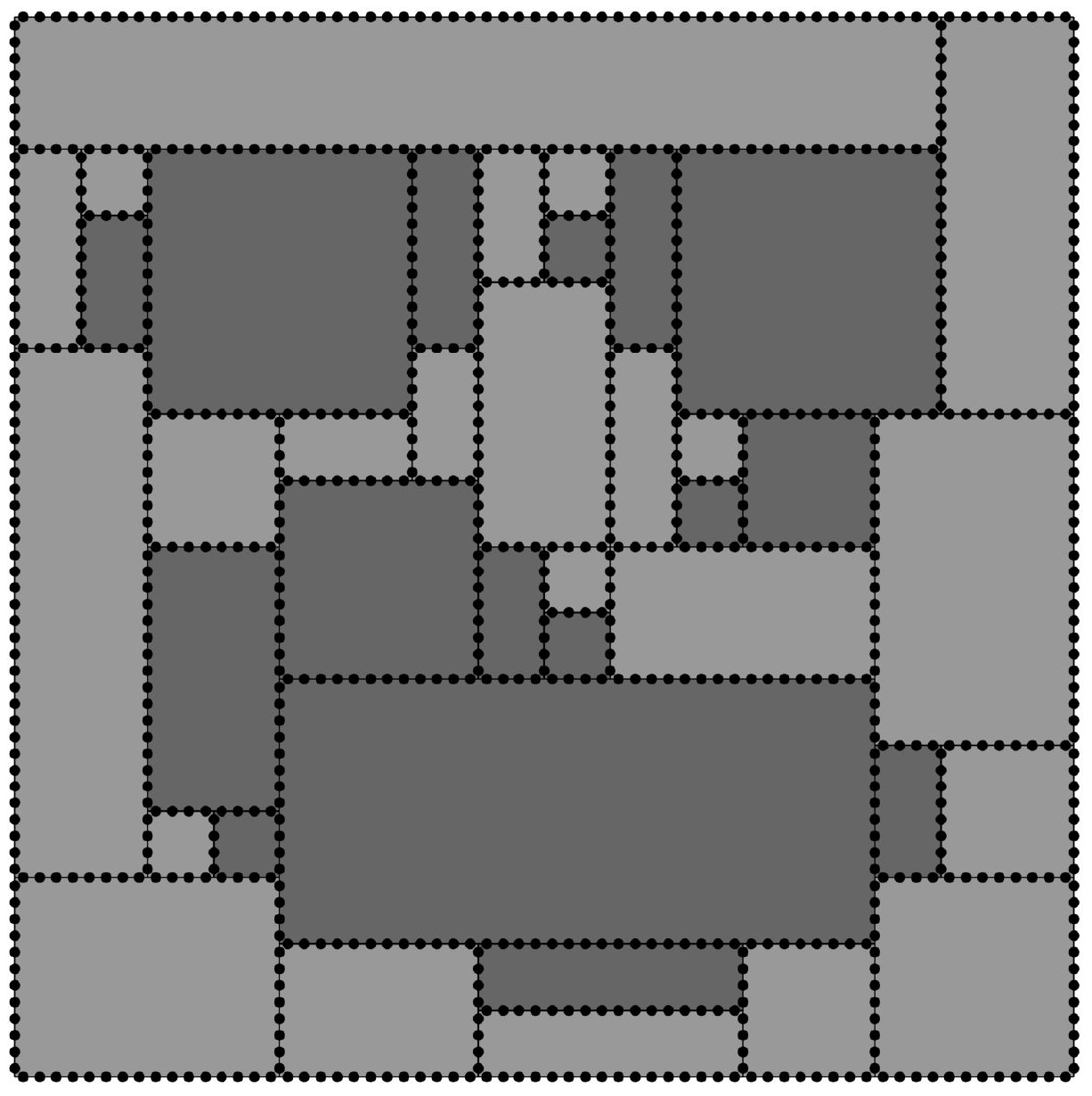}}
\caption{\revision{(a) The abscissas required for the semi-analytical method for the pixelated $16$ by $16$ geometry presented in Table \ref{tab:pixelation_refinement} with \amend{$N_x=N_y = 4$}, which would require the solution of a $2304$ by $2304$ linear \amend{system (\ref{eq:linear_system_paper2})} to compute the effective \amend{diffusivity}. (b) The abscissas for the extended semi-analytical method for the same geometry, in which rectangular blocks of different sizes are considered and interfaces do not have to be aligned. The spacing between abscissas is set to be $1/64$ and requires only the solution of a \amend{$803$} by \amend{$803$} linear system to compute the effective \amend{diffusivity}.}}
\label{fig:ESA}
\end{figure}

Although in this paper we consider geometries comprised of \revision{square blocks} in which all interfaces are aligned (such as in \amend{Figure \ref{fig:figure1}(b))}, our semi-analytical method can be extended \revision{to non-aligned rectangular blocks where \amend{interfaces} between adjacent blocks of homogeneous material are removed wherever possible. For example, in Figure \ref{fig:ESA}(b), we consider the pixelated geometry corresponding to $r = 16$ in Table \ref{tab:pixelation_refinement} and divide it into a number of rectangles of different sizes. The rectangles are chosen in such a way as to attempt to minimise the sum of the perimeter of the rectangles. The abscissas used for the numerical integration of the flux functions in equations (\ref{eq:quadrature})--(\ref{eq:quadrature2}) are placed on the interfaces between adjacent rectangles. Due to the fewer interfaces in Figure \ref{fig:ESA}(b), the total number of abscissas required is reduced, which would lead to decreased computational time, due to the smaller size of the linear system (\ref{eq:linear_system_paper2}) and increased accuracy, due to the fewer integrals along interfaces (\ref{eq:quadrature})--(\ref{eq:quadrature2}) required to be approximated.} For example in Figure \amend{\ref{fig:ESA}(a)}, we show the abscissas required for the semi-analytical method for the pixelated $16$ by $16$ geometry from Table \ref{tab:pixelation_refinement} with $N_x = N_y = 4$, which requires the solution of a $2304$ by $2304$ linear system in order to calculate the effective diffusivity. \amend{The} linear system corresponding to the formulation in \amend{Figure \ref{fig:ESA}(b), however is only }\amend{$803$ by $803$}, where \amend{$764$} of the unknowns in the linear system correspond to the abscissas of the flux function evaluations on the interfaces and the remaining \amend{$39$} unknowns correspond to the constants $K_{i,j}$ that appear in equation (\ref{eq:solution}). \amend{This \textit{extended semi-analytical method}} reduces the size of the linear system by over \amend{$65\%$} for the example in Figure \ref{fig:ESA} and highlights \amend{its potential} to greatly reduce the computational time required to calculate \edit{effective diffusivities.} \finaledit{\finalfinaledit{In this sense, the method could be considered as a form of adaptive approach}, as the abscissas used in the calculation of the effective diffusivities are chosen based on the diffusivities of the individual blocks and may change if the underlying geometry changes throughout the course of a simulation. The \textit{extended semi-analytical method} shares some of the advantages and disadvantages of numerical solution methods that use adaptive meshes that fit the geometry of the unit cell, such as the methods presented in \absolutefinaledit{\citep{sohn_2018} and \citep{amanbek_2019}.} All of these methods can be more computationally efficient due to fewer unknowns required to represent the solution $\psi^{(\xi)}$ in sections of the unit cell where the geometry is unchanged. However pre-processing is required to generate an adaptive mesh and remeshing may be required if the geometry changes throughout a simulation.}


\finalrevision{\subsection{Extension of semi-analytical method to discontinuous media}
\finaledit{Although in this work we consider \newedit{positive block diffusivities}, our method can be extended to discontinuous media in which some of the blocks have zero diffusivity.} To account for this, on any interface between two blocks in which at least one of the blocks has zero diffusivity, the interface and boundary functions (\ref{eq:BC1})--(\ref{eq:BC2}) would be set to equal zero along that interface. Additionally, no boundary value problem (\ref{eq:diffusion_2}) is imposed on any block in which the diffusivity is zero, as equation (\ref{eq:diffusion_2}) would be trivially satisfied. Any $m$ by $n$ geometry with blocks that have zero diffusivity would require fewer unknown interface functions than an $m$ by $n$ geometry without blocks that have zero diffusivity. This would decrease the size of the \newedit{linear system (\ref{eq:linear_system_paper2})}. For example, for the geometry presented in Figure \ref{fig:ESA}(a), with $N_x = N_y = 4$, the linear system has dimension $2304$ by $2304$. \newedit{However, if} the dark grey \neweredit{blocks had a} diffusivity of zero instead of $0.1$, the linear system \newedit{has} dimension $1155$ by $1155$. Furthermore, the extended semi-analytical method presented in section \ref{sec:esa} could \newedit{also} be applied to geometries which contain blocks with \newedit{zero diffusivity to further reduce the dimension of the linear system (\ref{eq:linear_system_paper2})}. For example, for the formulation presented in Figure \ref{fig:ESA}(b), if the dark grey blocks were set to have a diffusivity of zero, \newedit{removing the interfaces between dark grey and light grey blocks results in a linear system of dimension 327 by 327.}}
\begin{table}
\centering
 \begin{tabular}{|l|c|c|c|c|c|}
 
 \hline
 \multicolumn{5}{|c|}{Benchmark} \\
 \hline
 \multicolumn{5}{|c|}{} \\[-0.17cm]
  \multicolumn{5}{|c|}{{\includegraphics[height = 0.2\textwidth]{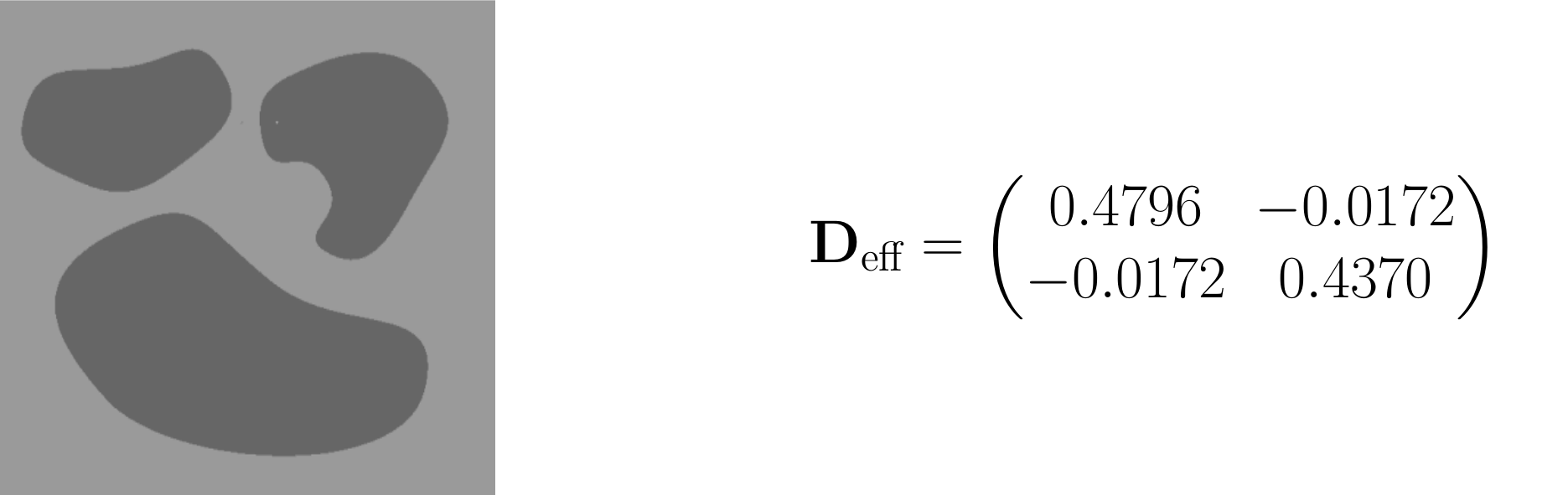}}} \\

\hline
Geometry & $\mathbf{D}_\text{eff}$  & $\mathbf{E}$ & \revision{$N$} & Time(s)\\
\hline
& & & & \\[-0.17cm]
{\includegraphics[width=0.18\textwidth,height = 0.18\textwidth]{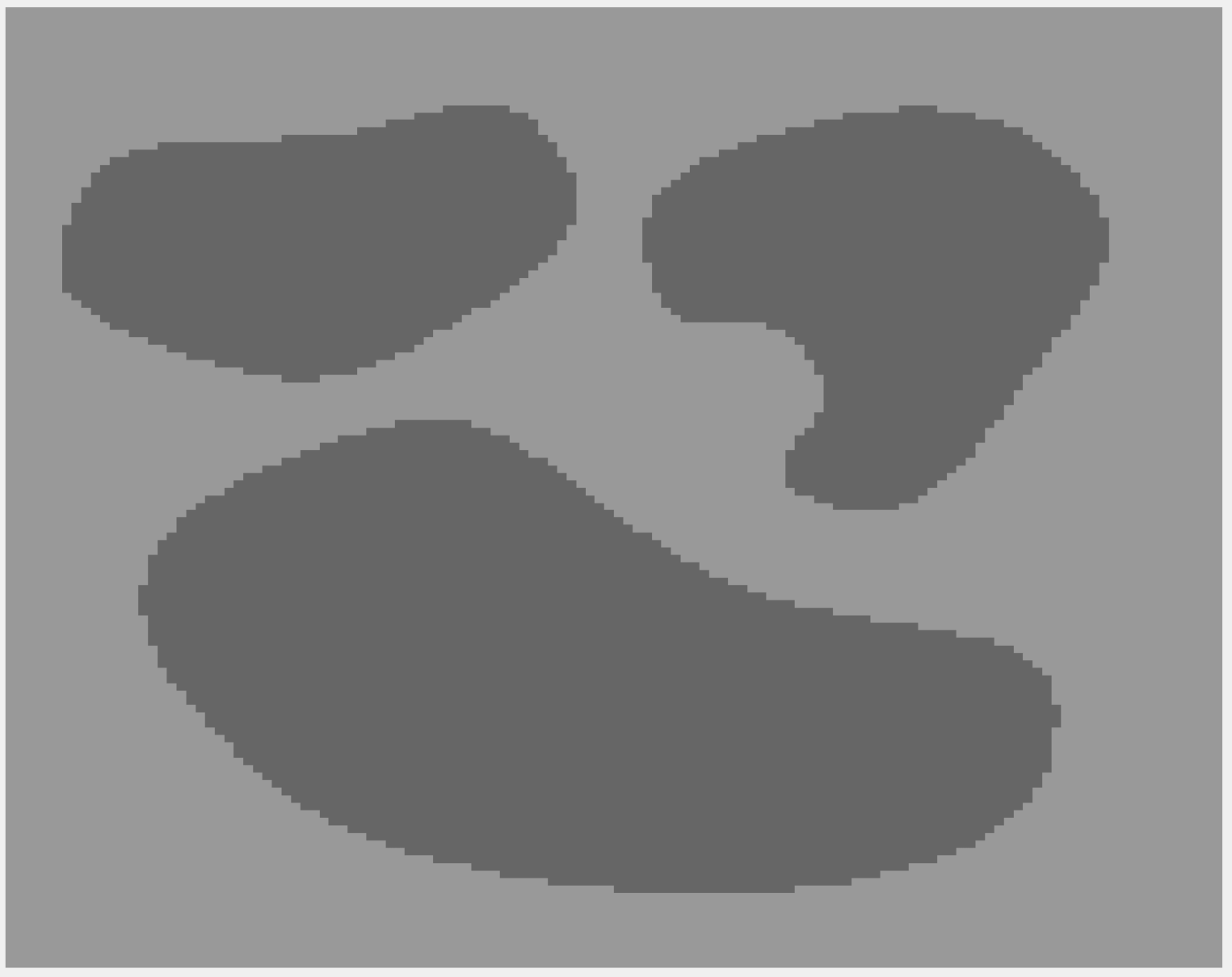}}  & \multirow{-8}{*}{{$ \footnotesize{\begin{pmatrix} 0.4761&  -0.0171 \\ -0.0171 & 0.4344\end{pmatrix}}$}} &  \multirow{-8}{*}{{$\footnotesize{\begin{pmatrix} 7.18 \text{e}{-03}&  6.12 \text{e}{-03}\\ 6.12 \text{e}{-03} & 5.79 \text{e}{-03}\end{pmatrix}} $}}&   \multirow{-8}{*}{$147456$} &  \multirow{-8}{*}{\edit{\amend{328}}}\\
$r = 128$ & & & &\\
\hline
& & & & \\[-0.17cm]
{\includegraphics[width=0.18\textwidth,height = 0.18\textwidth]{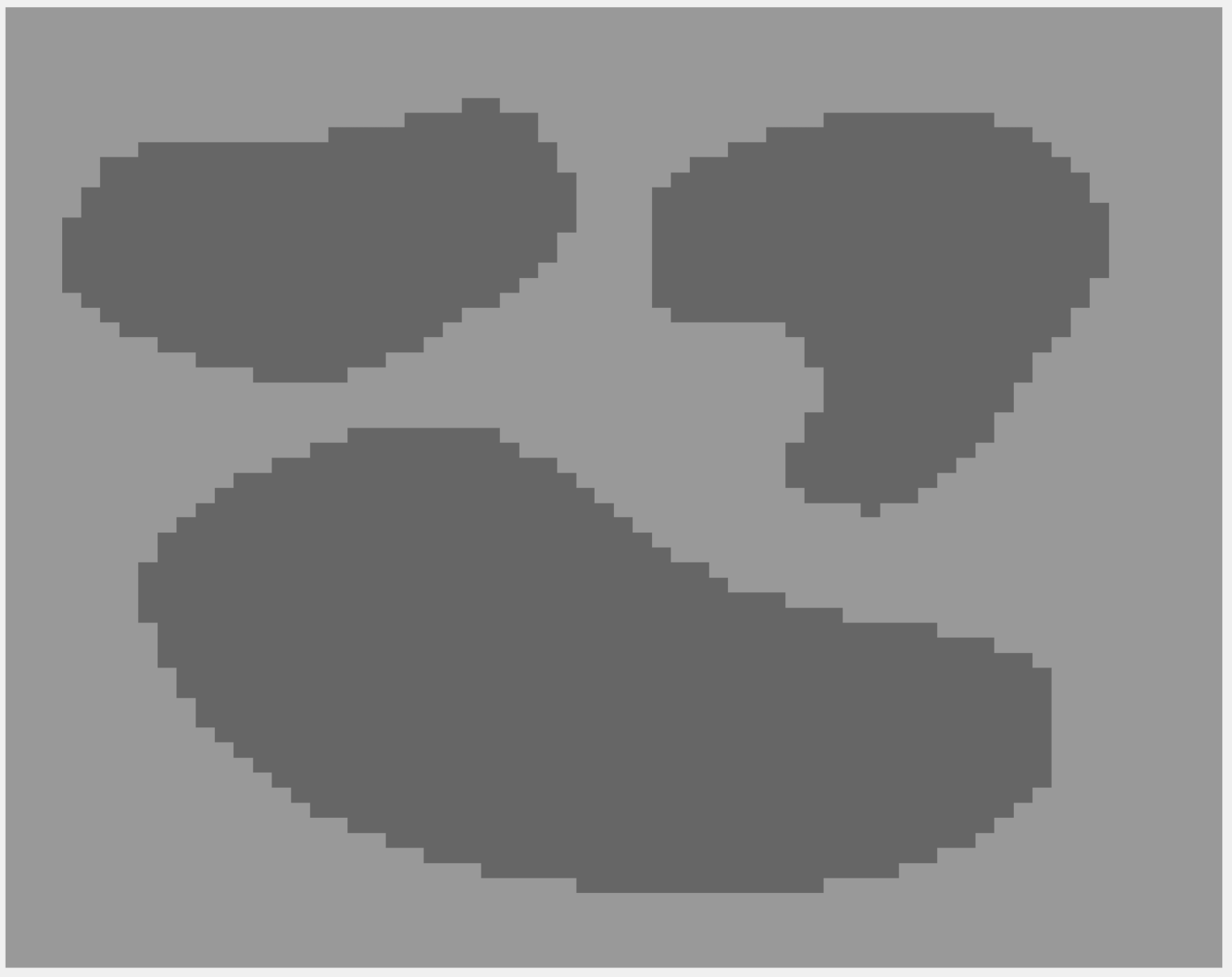}}&  \multirow{-8}{*}{{$\footnotesize{ \begin{pmatrix} 0.4740&  -0.0170 \\ -0.0170 & 0.4336\end{pmatrix}} $}} & \multirow{-8}{*}{{$ \footnotesize{\begin{pmatrix} 1.16 \text{e}{-02}&  1.44 \text{e}{-02}\\ 1.44 \text{e}{-02} & 7.81 \text{e}{-03}\end{pmatrix}} $}} &   \multirow{-8}{*}{$36864$} &  \multirow{-8}{*}{\edit{\amend{12.5}}}\\
$r = 64$ & & & &\\
\hline
& & & & \\[-0.17cm]
{\includegraphics[width=0.18\textwidth,height = 0.18\textwidth]{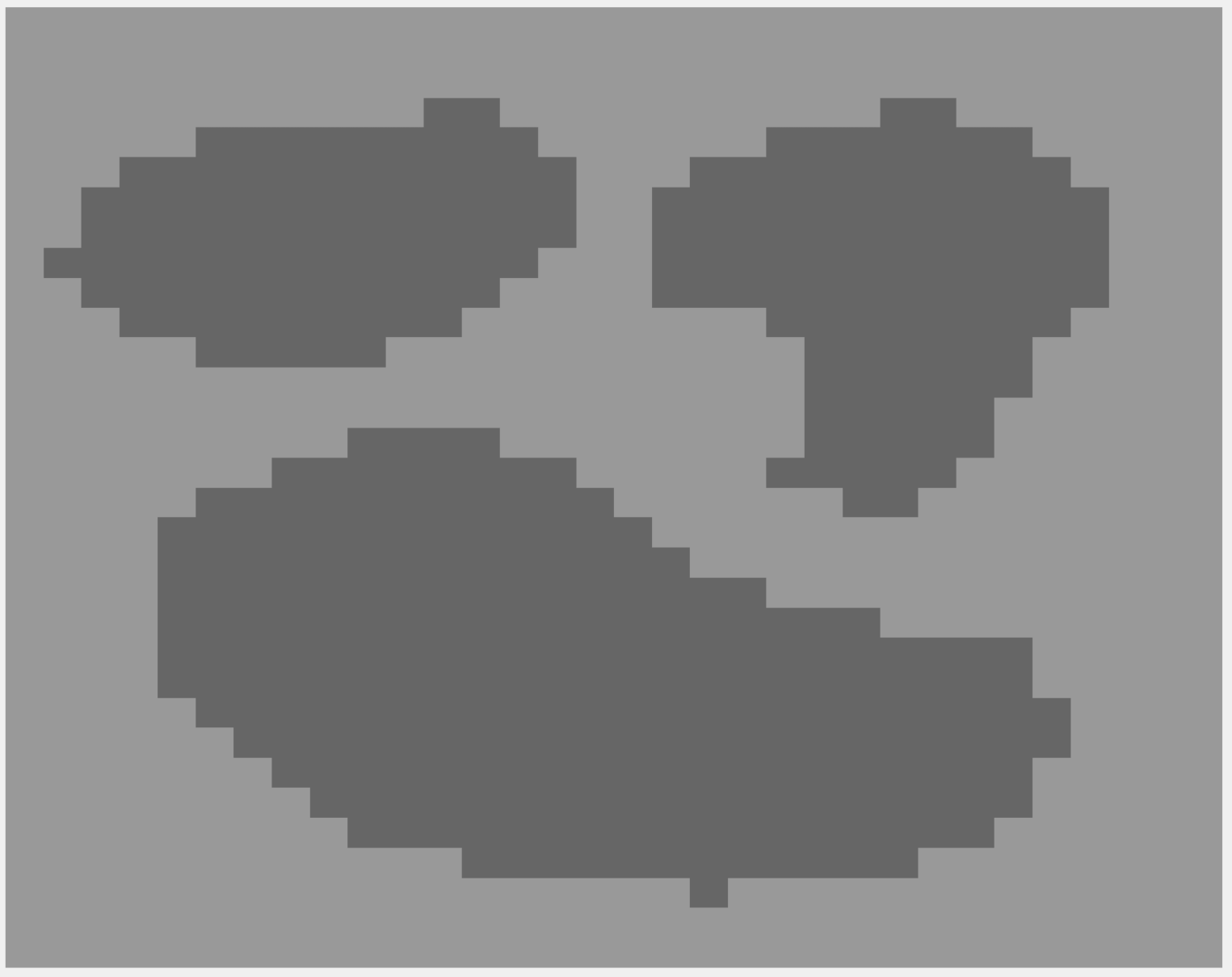}}&  \multirow{-8}{*}{{$\footnotesize{ \begin{pmatrix} 0.4750&  -0.0146 \\ -0.0146 & 0.4314\end{pmatrix}} $}} & \multirow{-8}{*}{{$ \footnotesize{\begin{pmatrix} 9.53 \text{e}{-03}&  1.49 \text{e}{-01}\\ 1.49 \text{e}{-01} & 1.27 \text{e}{-02}\end{pmatrix}} $}} &  \multirow{-8}{*}{$9216$}&  \multirow{-8}{*}{\amend{0.912}}\\ 
$r = 32$ & & & &\\
\hline
& & & & \\[-0.17cm]
{\includegraphics[width=0.18\textwidth,height = 0.18\textwidth]{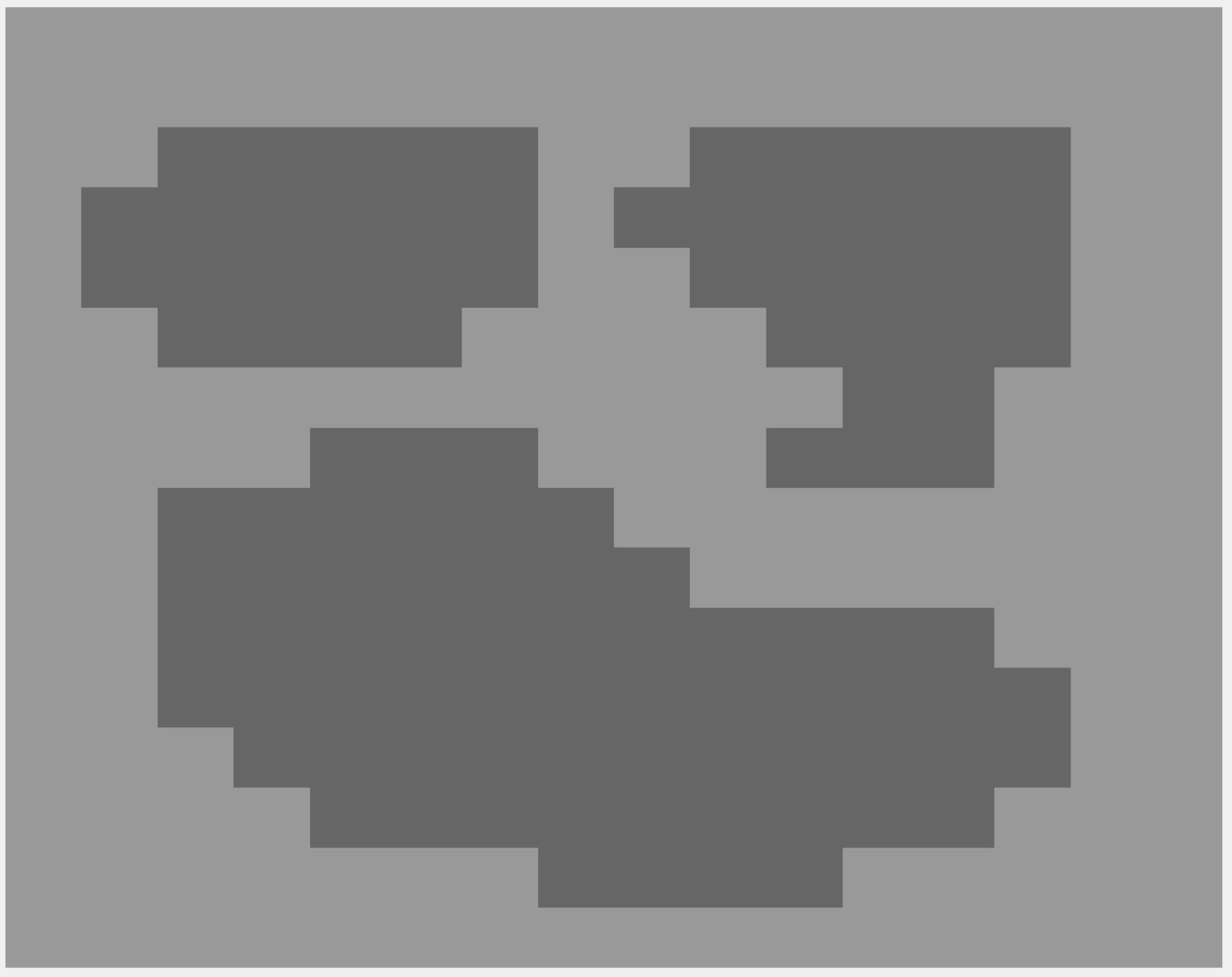}}&  \multirow{-8}{*}{{$\footnotesize{ \begin{pmatrix} 0.4689&  -0.0194 \\ -0.0194 & 0.4362\end{pmatrix}} $}} &  \multirow{-8}{*}{{$\footnotesize{ \begin{pmatrix} 2.22 \text{e}{-02}&  1.29 \text{e}{-01}\\ 1.29 \text{e}{-01} & 1.66 \text{e}{-03}\end{pmatrix}} $}} &   \multirow{-8}{*}{$2304$}&  \multirow{-8}{*}{\amend{0.0852}}\\ 
$r = 16$ & & & &\\
\hline
\end{tabular}
\caption{\revision{Pixelated geometries corresponding to \amend{Figure \ref{fig:figure1}(a)} with corresponding effective \edit{diffusivities} ($\mathbf{D}_\text{eff}$), relative errors ($\mathbf{E}$), number of rows of blocks in the pixelated geometry($r$), \amend{number of rows in the corresponding linear system (\ref{eq:linear_system_paper2}) ($N$)} and time taken to solve the linear system (Time\edit{(s)}).}}
\label{tab:pixelation_refinement}
\end{table}

\section{Conclusion}
\label{sec:conclusion_paper2}
\finaledit{In this \final{paper} we have presented \revision{a} semi-analytical method {for solving the classical boundary value problem arising from the homogenization of a} two-dimensional, {pixelated, locally-isotropic, heterogeneous, periodic domain}. \revision{Our} approach {involves reformulating} the heterogeneous boundary value problem on the unit cell as a family of {boundary value problems on the homogeneous blocks} by introducing unknown functions representing the flux at the interfaces between adjacent blocks. We then solve each homogeneous boundary value problem via standard techniques to yield an analytical solution \newedit{depending on the integrals involving} unknown interface functions. By applying an appropriate numerical quadrature rule to these integrals and enforcing continuity of the solution across each interface, {the solution of the original heterogeneous boundary value problem can be computed, which allows the effective diffusivity to be calculated}. }{Our numerical experiments demonstrated \amend{our} semi-analytical method yields accurate results when applied} to standard test problems from the literature. \finaledit{For all of these test problems\finalfinaledit{ considered in this work}, we found that our semi-analytical method is faster and/or more accurate than {a standard} finite volume method.} {The primary reason for this is that the linear system for the semi-analytical method is much smaller than the linear system for the finite volume method as unknowns are located only along interfaces between adjacent blocks and not in the interior of the blocks \revision{(see Figure \ref{fig:figure1}(d)--(e))}. We also} demonstrated that our semi-analytical method can be applied to complex heterogeneous \newedit{geometries by pixellating and representing the geometry as a grid of blocks, as demonstrated in Table \ref{tab:pixelation_refinement}.}

While all of the heterogeneous media considered in this paper consist of only two different materials, our semi-analytical method is not restricted to such media \revision{as the diffusivity in each block can be different to all other blocks}. {Furthermore, we only considered periodic boundary conditions (\ref{eq:UBC5})--(\ref{eq:UBC8}), however, our semi-analytical method can be modified to accomodate additional forms of boundary conditions used in other homogenization \amend{techniques \citep{szymkiewicz_2012} and other }steady-state diffusion problems} on two-dimensional, block heterogeneous media. {There is also} room for improvement in our semi-analytical method that may further {improve} its efficiency and accuracy. For example: (i) allowing for blocks of different sizes where not all interfaces are {aligned would} reduce the number of interface functions and thus reduce the size of the linear system (\ref{eq:linear_system_paper2}), as discussed in \amend{section \ref{sec:esa}}; (ii) applying a quadrature rule to approximate the integrals involving the interface functions that allows for more than $2N_x-3$ terms to be taken in the \revision{summations} {would likely improve the accuracy}{; (iii) applying an iterative technique to the solution of the linear system, as demonstrated in \citep{vondrejc_2014}, may allow for improvement in the efficiency of our semi-analytical method.} \finaledit{In \finalfinaledit{our future work}, we plan to investigate the effect of the accuracy of an effective diffusivity tensor on the accuracy of the coarse-scale solution of diffusion problems in complex heterogeneous media and the use of our semi-analytical method for alleviating bottlenecks in coarse-grained models of non-periodic heterogeneous media, where the domain is decomposed into several sub-domains, which are individually homogenized.}



\section*{Acknowledgments}
The first author acknowledges discussions with Nadja Ray and Andreas Rupp from Friedrich-Alexander University of Erlangen-N{\"u}rnberg on randomly generating aggregated heterogeneous domains. \finaledit{All authors acknowledge the helpful comments of the anonymous reviewers and editors that helped improve the quality of the manuscript}. The second and third authors acknowledge funding from the Australian Research Council (DE150101137, DP150103675).

\appendix

\section{Boundary conditions for blocks in bottom row and right column}
\label{sec:BCs}
{The boundary conditions (\ref{eq:BC1})--(\ref{eq:BC2}) are valid for all blocks except those in either the bottom row ($i = m$) or right column $(j = n)$. For the blocks in the bottom row and right column, the boundary conditions are}:
\begin{gather}
\label{eq:BC3}
D_{m,j}\frac{\partial v_{m,j}^{(\xi)}}{\partial x}(x_{j-1},y) = g_{(m-1)n+j}(y), \quad  y_{m-1} < y < y_m,\\
D_{m,j}\frac{\partial v_{m,j}^{(\xi)}}{\partial x}(x_j,y) = g_{(m-1)n+j+1}(y), \quad y_{m-1} < y < y_m,\\
\label{eq:BC4}
D_{m,j}\frac{\partial v_{m,j}^{(\xi)}}{\partial y}(x,y_{m-1}) = q_{jm}(x), \quad x_{j-1} < x < x_j,\\
D_{m,j}\frac{\partial v_{m,j}^{(\xi)}}{\partial y}(x,y_m) = q_{(j-1)m+1}(x), \quad x_{j-1} < x < x_j,
\end{gather}
for $j = 1,\hdots,n-1$,
\begin{gather}
\label{eq:BC5}
D_{i,n}\frac{\partial v_{i,n}^{(\xi)}}{\partial x}(x_{n-1},y) = g_{in}(y), \quad y_{i-1} < y < y_i,\\
D_{i,n}\frac{\partial v_{i,n}^{(\xi)}}{\partial x}(x_n,y) = g_{(i-1)n+1}(y), \quad y_{i-1} < y < y_i,\\
\label{eq:BC6}
D_{i,n}\frac{\partial v_{i,n}^{(\xi)}}{\partial y}(x,y_{i-1}) = q_{(n-1)m+i}(x),  \quad x_{n-1} < x < x_n,\\
D_{i,n}\frac{\partial v_{i,n}^{(\xi)}}{\partial y}(x,y_i) = q_{(n-1)m+i+1}(x), \quad x_{n-1} < x < x_n,
\end{gather}
for $i = 1,\hdots,m-1$, {and}
\begin{gather}
\label{eq:BC7}
D_{m,n}\frac{\partial v_{m,n}^{(\xi)}}{\partial x}(x_{n-1},y) = g_{mn}(y), \quad y_{m-1} < y < y_m,\\
D_{m,n}\frac{\partial v_{m,n}^{(\xi)}}{\partial x}(x_n,y) = g_{(m-1)n+1}(y), \quad y_{m-1} < y < y_m,\\
\label{eq:BC8}
D_{m.n}\frac{\partial v_{m,n}^{(\xi)}}{\partial y}(x,y_{m-1}) =q_{mn}(x), \quad x_{n-1} < x < x_n,\\ 
D_{m.n}\frac{\partial v_{m,n}^{(\xi)}}{\partial y}(x,y_m) =q_{(n-1)m+1}(x), \quad x_{n-1} < x < x_n.
\end{gather}
\section{Coefficients for blocks in bottom row and right column}
\label{sec:coefficients}
{The coefficients (\ref{eq:coeff_1_paper2})--(\ref{eq:coeff_4_paper2}) are valid for all blocks except those in either the bottom row ($i = m$) or right column $(j = n)$.  For the blocks in the bottom row and right column, the coefficients are defined as}:
\begin{gather}
a_{m,j,k}^{(\xi)} = \frac{2}{h_m} \int_{y_{m-1}}^{y_m} \frac{g_{(m-1)n+j}(y)}{D_{m,j}} \cos\left(\frac{k\pi(y-y_{m-1})}{h_m}\right)\, dy, \\
b_{m,j,k}^{(\xi)} = \frac{2}{h_m} \int_{y_{m-1}}^{y_m} \frac{g_{(m-1)n+j+1}(y)}{D_{m,j}} \cos\left(\frac{k\pi (y-y_{m-1})}{h_m}\right)\, dy,\\
c_{m,j,k}^{(\xi)} = \frac{2}{l_j} \int_{x_{j-1}}^{x_j} \frac{q_{jm}(x)}{D_{m,j}} \cos\left(\frac{k\pi (x-x_{j-1})}{l_j}\right)\, dx,\\
d_{m,j,k}^{(\xi)} = \frac{2}{l_j} \int_{x_{j-1}}^{x_j} \frac{q_{(j-1)m+1}(x)}{D_{m,j}} \cos\left(\frac{k\pi (x-x_{j-1})}{l_j}\right)\, dx,
\end{gather}
for $j = 1,\hdots,n-1$,
\begin{gather}
a_{i,n,k}^{(\xi)} = \frac{2}{h_i} \int_{y_{i-1}}^{y_i} \frac{g_{in}(y)}{D_{i,n}} \cos\left(\frac{k\pi (y-y_{i-1})}{h_i}\right)\, dy,\\
b_{i,n,k}^{(\xi)} = \frac{2}{h_i} \int_{y_{i-1}}^{y_i} \frac{g_{(i-1)n+1}(y)}{D_{i,n}} \cos\left(\frac{k\pi (y-y_{i-1})}{h_i}\right)\, dy,\\
c_{i,n,k}^{(\xi)} = \frac{2}{l_n} \int_{x_{n-1}}^{x_n} \frac{q_{(n-1)m+i}(x)}{D_{i,n}} \cos\left(\frac{k\pi(x-x_{n-1})}{l_n}\right)\, dx,\\
d_{i,n,k}^{(\xi)} = \frac{2}{l_n} \int_{x_{n-1}}^{x_n} \frac{q_{(n-1)m+i+1}(x)}{D_{i,n}} \cos\left(\frac{k\pi (x-x_{n-1})}{l_n}\right)\, dx,
\end{gather}
for $i = 1,\hdots,m-1$, {and}
\begin{gather}
a_{m,n,k}^{(\xi)} = \frac{2}{h_m} \int_{y_{i-1}}^{y_i} \frac{g_{mn}(y)}{D_{m,n}} \cos\left(\frac{k\pi (y-y_{m-1})}{h_m}\right)\, dy,\\
b_{m,n,k}^{(\xi)} = \frac{2}{h_m} \int_{y_{i-1}}^{y_i} \frac{g_{(m-1)n+1}(y)}{D_{m,n}} \cos\left(\frac{k\pi (y-y_{m-1})}{h_m}\right)\, dy,\\
c_{m,n,k}^{(\xi)} = \frac{2}{l_n} \int_{x_{n-1}}^{x_n} \frac{q_{mn}(x)}{D_{m,n}} \cos\left(\frac{k\pi(x-x_{n-1})}{l_n}\right)\, dx,\\
d_{m,n,k}^{(\xi)} = \frac{2}{l_n} \int_{x_{n-1}}^{x_n} \frac{q_{(n-1)m+1}(x)}{D_{m,n}} \cos\left(\frac{k\pi (x-x_{n-1})}{l_n}\right)\, dx.
\end{gather}
\section{Solvability conditions for blocks in bottom row and right column}
\label{sec:solvability}
{The solvability conditions (\ref{eq:S1}) are valid for all blocks except those in either the bottom row ($i = m$) or right column $(j = n)$. For the blocks in the bottom row and right column, the solvability conditions are:}
\begin{multline}
\label{eq:S2}
\int_{y_{m-1}}^{y_m} \frac{g_{(m-1)n+j}(y)}{D_{m,j}}\, dy - \int_{y_{m-1}}^{y_m} \frac{g_{(m-1)n+j+1}(y)}{D_{m,j}}\, dy \\+  \int_{x_{j-1}}^{x_j} \frac{q_{jm}(x)}{D_{m,j}} \, dx - \int_{x_{j-1}}^{x_j} \frac{q_{(j-1)m+1}(x)}{D_{m,j}} \, dx = 0,
\end{multline}
for $j = 1,\hdots,n-1$,
\begin{multline}
\label{eq:S3}
\int_{y_{i-1}}^{y_i} \frac{g_{in}(y)}{D_{i,n}}\, dy - \int_{y_{i-1}}^{y_i} \frac{g_{(i-1)n+1}(y)}{D_{i,n}}\, dy \\+  \int_{x_{n-1}}^{x_n} \frac{q_{(n-1)m+i}(x)}{D_{i,n}} \, dx - \int_{x_{n-1}}^{x_n} \frac{q_{(n-1)m+i+1}(x)}{D_{i,n}} \, dx = 0,
\end{multline}
for $i = 1,\hdots,m-1$, {and}
\begin{multline}
\label{eq:S4}
\int_{y_{m-1}}^{y_m} \frac{g_{mn}(y)}{D_{m,n}}\, dy - \int_{y_{m-1}}^{y_m} \frac{g_{(m-1)n+1}(y)}{D_{m,n}}\, dy \\+  \int_{x_{n-1}}^{x_n} \frac{q_{mn}(x)}{D_{m,n}} \, dx - \int_{x_{n-1}}^{x_n} \frac{q_{(n-1)m+1}(x)}{D_{m,n}} \, dx = 0.
\end{multline}

\raggedright

\bibliographystyle{elsarticle-num}
\bibliography{mypapers}

\end{document}